\PassOptionsToPackage{inline}{enumitem}

\documentclass[submission]{Arxiv}

\usepackage[english]{babel}
\usepackage{amsmath,amsfonts,amssymb,amsthm,bm}
\usepackage{mathtools}
\usepackage[T1]{fontenc}
\usepackage{dsfont}
\usepackage{multicol}
\usepackage{stmaryrd}
\usepackage{mleftright}
\usepackage{tikz}
\usetikzlibrary{arrows.meta,fit}

\usepackage{xcolor}
\definecolor{ColA}{HTML}{0061b7} 
\definecolor{ColB}{HTML}{b75600} 
\newcommand{\ColMix}[2]{\textcolor{ColA!#1!ColB}{#2}}

\setenumerate{noitemsep,topsep=3pt}

\usepackage[backend=bibtex]{biblatex}
\numberwithin{equation}{section}

\newtheorem{Theorem}{Theorem}[section]
\newtheorem{Proposition}[Theorem]{Proposition}

\renewcommand{\leq}{\leqslant}
\renewcommand{\geq}{\geqslant}

\newcommand{\Hide}[1]{\ColMix{100}{\tt HIDEN}}
\newcommand{\Def}[1]{\ColMix{20}{\em #1}}
\newcommand{\Par}[1]{\mleft(#1\mright)}
\newcommand{\Bra}[1]{\mleft\{#1\mright\}}
\newcommand{\Han}[1]{\mleft[#1\mright]}

\newcommand{\Angle}[1]{\mleft\langle#1\mright\rangle}

\newcommand{\BBrack}[1]{\bm{[}#1\bm{]}}

\newcommand{\Bag}[1]{\lbag#1\rbag}

\tikzstyle{Centering}=[{baseline={([yshift=-0.5ex]current bounding box.center)}}]

\title{Polynomial realizations of natural Hopf algebras of nonsymmetric operads}
\author{%
    Samuele Giraudo%
    \thanks{
        \href{mailto:giraudo.samuele@uqam.ca}{\tt giraudo.samuele@uqam.ca}.
        This research has been partially supported by the projects CARPLO (ANR-20-CE40-0007)
        and ALCOHOL (ANR-19-CE40-0006) of the Agence nationale de la recherche.}%
    \addressmark{1}}
\address{%
    \addressmark{1} Université du Québec à Montréal, LACIM, Pavillon Président-Kennedy, 201
    Avenue du Président-Kennedy, Montréal, H2X~3Y7, Canada.}
\abstract{%
    The natural Hopf algebra $\NaturalHopfAlgebra \App \Operad$ of an operad $\Operad$ is a
    Hopf algebra whose bases are indexed by some words on $\Operad$. We introduce new bases
    of these Hopf algebras deriving from free operads via new lattice structures on their
    basis elements. We construct polynomial realizations of $\NaturalHopfAlgebra \App
    \Operad$ by using alphabets of variables endowed with unary and binary relations. By
    specializing our polynomial realizations, we discover links between $\NaturalHopfAlgebra
    \App \Operad$ and the Hopf algebra of word quasi-symmetric functions of Hivert, the
    decorated versions of the Connes-Kreimer Hopf algebra of Foissy, the Faà di Bruno Hopf
    algebra and its deformations, and the noncommutative multi-symmetric functions Hopf
    algebras of Novelli and Thibon.}
\received{\today}

\newcommand{\K}{\mathbb{K}}
\newcommand{\N}{\mathbb{N}}
\newcommand{\GenA}{\mathsf{a}}
\newcommand{\GenB}{\mathsf{b}}
\newcommand{\GenC}{\mathsf{c}}
\newcommand{\GenG}{\mathsf{g}}
\newcommand{\Iverson}[1]{\BBrack{#1}}
\newcommand{\App}{\mathop{\cdot}}
\newcommand{\Length}{\ell}

\newcommand{\Degree}{\mathrm{dg}}
\newcommand{\Arity}{\mathrm{ar}}
\newcommand{\Conc}{\mathbin{.}}
\newcommand{\Leaf}{\perp}
\newcommand{\TermT}{\mathfrak{t}}
\newcommand{\TermS}{\mathfrak{s}}
\newcommand{\ForestF}{\mathfrak{f}}
\newcommand{\Signature}{\mathcal{S}}
\newcommand{\SignatureExample}{\Signature_{\mathrm{e}}}
\newcommand{\SetTerms}{\mathfrak{T}}
\newcommand{\SetForests}{\mathfrak{F}}
\newcommand{\Operad}{\mathcal{O}}
\newcommand{\Corolla}{\iota}
\newcommand{\As}{\mathsf{As}}
\newcommand{\Unit}{\mathds{1}}
\newcommand{\Reduced}{\mathrm{rd}}
\newcommand{\Edge}[2]{\mathrel{\leadsto}_{#1, #2}}
\newcommand{\EasterlyWindCovering}{\rightharpoonup}
\newcommand{\Leq}{\mathrel{\preccurlyeq}}
\newcommand{\Equiv}{\mathrel{\equiv}}
\newcommand{\EasterlyWindPoset}{\mathcal{W}}
\newcommand{\Mobius}{\mu}
\newcommand{\Product}{\mathbin{\cdot}}
\newcommand{\DummyConstant}{\mathfrak{d}}
\newcommand{\Space}{\mathcal{V}}
\newcommand{\NaturalHopfAlgebra}{\mathbf{N}}
\newcommand{\BasisE}{\mathsf{E}}
\newcommand{\BasisF}{\mathsf{F}}
\newcommand{\BasisH}{\mathsf{H}}
\newcommand{\BasisM}{\mathsf{M}}
\newcommand{\DecorationRelation}{\mathrm{D}}
\newcommand{\RootRelation}{\mathrm{R}}
\newcommand{\EdgeRelation}[1]{\mathrel{\prec_{#1}}}
\newcommand{\AlphabetSum}{\mathbin{+\mkern-13mu+}}
\newcommand{\CompatibleWord}{\mathrel{\Vdash}}
\newcommand{\CompatiblePackedWord}{\mathrel{\vdash}}
\newcommand{\Realization}{\mathsf{r}}
\newcommand{\Coproduct}{\Delta}
\newcommand{\WordToTensor}{\theta}
\newcommand{\Alphabet}{\mathbb{A}}
\newcommand{\AlphabetN}{\Alphabet_{\N}}
\newcommand{\Letter}{\mathbf{a}}

\newcommand{\WQSym}{\mathbf{WQSym}}
\newcommand{\NCK}{\mathbf{NCK}}
\newcommand{\NCFdB}{\mathbf{FdB}}
\newcommand{\NCSF}{\mathbf{Sym}}
\newcommand{\Pack}{\mathrm{pck}}
\newcommand{\SetPackedWords}{\mathcal{P}}

\newcommand{\MAs}{\mathsf{MAs}}
\newcommand{\MultisetM}{\mathfrak{m}}

\newcommand{\NAryRelation}{\mathcal{R}}
\newcommand{\HopfAlgebra}{\mathcal{H}}
\newcommand{\ClassAlphabets}{\mathcal{A}}

\DeclareMathOperator{\Over}{%
    \begin{tikzpicture}[Centering,scale=.16]
        \draw(0,0)--(1,1.25);
    \end{tikzpicture}}

\DeclareMathOperator{\Under}{%
    \begin{tikzpicture}[Centering,scale=.16]
        \draw(0,0)--(1,-1.25);
    \end{tikzpicture}}

\tikzstyle{Node}=[circle,draw=ColA!90,fill=ColA!5,inner sep=1pt,minimum size=3.5mm,thick,
font=\scriptsize]
\tikzstyle{Leaf}=[circle,draw=ColA!90,fill=ColA!5,inner sep=0pt,minimum size=1mm,thick]
\tikzstyle{Edge}=[draw=ColB!95,cap=round,rounded corners=2.5pt,thick]
\tikzstyle{FitNodes}=[draw=ColB!100,fill=ColB!40,opacity=.25,inner sep=2pt,very
thick,rounded corners=6pt]
\tikzstyle{GraphEdge}=[ColB,cap=round,thick]
\tikzstyle{GraphLabeledVertex}=[font=\footnotesize,node distance=3mm]

\begin{document}
\maketitle

\section*{Introduction}
A polynomial realization (abbreviated as \textit{PR}) of a Hopf algebra consists of
interpreting its elements as polynomials, either commutative or not, in such a way that its
product translates as polynomial multiplication and the coproduct translates as a simple
transformation of the alphabet of variables. A great portion of combinatorial Hopf algebras
appearing in combinatorics admit PRs~\cite{DHT02,NT06,FNT14,Foi20}. It is striking to note
that Hopf algebras involving a variety of different families of combinatorial objects and
operations on them can be translated and understood in a common manner through adequate PRs.

Such realizations are crucial for several reasons. First, they allow us to prove easily that
the Hopf algebra axioms (like associativity, coassociativity, and the Hopf compatibility
between the product and the coproduct) are satisfied. Indeed, if a Hopf algebra admits a PR,
such properties are almost immediate on polynomials~\cite{Hiv07}. Second, given a PR of a
Hopf algebra, it is in most cases fruitful to specialize the associated polynomials (for
instance by letting the variables commute) in order to get Hopf algebra morphisms to other
Hopf algebras. This leads to the construction of new Hopf algebras or establishes links
between already existing ones~\cite{FNT14}. Finally, such families of polynomials realizing
a Hopf algebra lead to the definition of new classes of special functions, analogous to
Schur or Macdonald functions~\cite{DHT02}.

Under right conditions, an operad $\Operad$ produces a Hopf algebra $\NaturalHopfAlgebra
\App \Operad$, called the \textit{natural Hopf algebra} of $\Operad$. The bases of
$\NaturalHopfAlgebra \App \Operad$ are indexed by some words on $\Operad$, and its coproduct
is inherited from the composition map of $\Operad$. This construction is considered
in~\cite{vdl04,CL07,BG16}, and we focus here on a noncommutative variation for nonsymmetric
operads~\cite{ML14}. In contrast to many examples of Hopf algebras having PRs, none are
known for $\NaturalHopfAlgebra \App \Operad$. The main contribution of this work is to
provide a PR for this family of Hopf algebras built from nonsymmetric operads. The
particularity of our approach is that we consider a PR based on variables belonging to
alphabets endowed with several unary and binary relations in order to capture the
particularities of the coproduct of $\NaturalHopfAlgebra \App \Operad$. This approach, using
what we call \textit{related alphabets}, generalizes the previous approaches using totally
ordered alphabets~\cite{DHT02,NT06,Hiv07}, quasi-ordered alphabets~\cite{Foi20}, or
alphabets endowed with a single binary relation~\cite{FNT14,Gir11}.

This work is presented as follows. In Section~\ref{sec:preliminaries}, the main notions
about natural Hopf algebras of operads, PRs, and free operads are provided. We introduce in
Section~\ref{sec:alternative_bases} two new bases of $\NaturalHopfAlgebra \App \Operad$ in
the case where $\Operad$ is a free operad. These constructions use new lattices on treelike
structures. In Section~\ref{sec:hopf_subalgebras}, we describe an injection of
$\NaturalHopfAlgebra \App \Operad /_{\Equiv}$ into $\NaturalHopfAlgebra \App \Operad$ when
$\Operad$ is a free operad and $\Equiv$ is an operad congruence subjected to some
properties. Section~\ref{sec:polynomial_realization} contains a PR of the natural Hopf
algebra of a free operad $\Operad$ and of some of its quotients. In the subsequent
Sections~\ref{sec:wqsym}, \ref{sec:connes_kreimer}, \ref{sec:faa_di_bruno},
and~\ref{sec:multi_symmetric_functions}, we present some consequences of our PR, leading to
links between $\NaturalHopfAlgebra \App \Operad$ and other Hopf algebras. In particular, we
show that a specialization of the PR of $\NaturalHopfAlgebra \App \Operad$ forms a Hopf
subalgebra of $\WQSym$, the Hopf algebra word quasi-symmetric functions~\cite{Hiv99,NT06}.
We provide, by building a free operad $\Operad_D$ depending on a set $D$ of decorations, a
construction of the Hopf algebra of $D$-decorated forests~\cite{Foi02a} as a specialization
of the PR of $\NaturalHopfAlgebra \App \Operad_D$. We show that the noncommutative Faà di
Bruno $\NCFdB$ Hopf algebra and its one-parameter deformation $\NCFdB_r$~\cite{Foi08} can be
built from the multi-multiassociative operads $\MAs_\Signature$, extending a construction
of~\cite{Gir16}. Finally, we build from $\MAs_\Signature$ a Hopf algebra $\NCFdB_{r,s}$
depending on two parameters such that $\NCFdB_{0, s}$ is the noncommutative multi-symmetric
function Hopf algebra~\cite{NT10} and $\NCFdB_{r, 0}$ is $\NCFdB_r$. By using our framework,
we obtain PRs of these Hopf algebras.

{\it General notations and conventions.}
If $f$ is an entity parameterized by an input $x$, we write $f \App x$ for $f(x)$. As $\App$
associates from right to left, $f \App g \App x$ denotes $f \App (g \App x)$. For any
statement $P$, the Iverson bracket $\Iverson{P}$ takes $1$ as value if $P$ is true and $0$
otherwise. For any integer $i$, $[i]$ denotes the set $\{1, \dots, i\}$. For any set $A$,
$A^*$ is the set of words on $A$. For any $w \in A^*$, $\Length \App w$ is the length of
$w$, and for any $i \in [\Length \App w]$, $w \App i$ is the $i$-th letter of $w$. The only
word of length $0$ is the empty word $\epsilon$. For any subset $A'$ of $A$, $w_{|A'}$ is
the subword of $w$ made of letters of $A'$. Given two words $w$ and $w'$, the concatenation
of $w$ and $w'$ is denoted by $ww'$ or by $w \Conc w'$.

\section{Natural Hopf algebras from operads} \label{sec:preliminaries}

\subsection{Natural Hopf algebras from nonsymmetric operads}
\label{subsec:natural_hopf_algebras}
We follow the usual notations about nonsymmetric operads~\cite{Gir18} (called simply
\Def{operads} here). Let $\Operad$ be an operad with arity map $\Arity$, composition map
$-\Han{-, \dots, -}$, and unit~$\Unit$. When each $x \in \Operad$ admits finitely many
factorizations $x = y \Han{z_1, \dots, z_{\Arity \App y}}$ where $y, z_1, \dots, z_{\Arity
\App y} \in \Operad$, $\Operad$ is \Def{finitely factorizable}. When there exists a map
$\Degree : \Operad \to \N$ such that $\Degree^{-1} \App 0 = \Bra{\Unit}$ and, for any $y,
z_1, \dots, z_{\Arity \App y} \in \Operad$, $\Degree \App y \Han{z_1, \dots, z_{\Arity \App
y}} = \Degree \App y + \Degree \App z_1 + \dots + \Degree \App z_{\Arity \App y}$, the map
$\Degree$ is a \Def{grading} of $\Operad$.

The \Def{natural Hopf algebra}~\cite{vdl04,CL07,ML14,BG16} of a finitely factorizable operad
$\Operad$ admitting a grading $\Degree$ is the Hopf algebra $\NaturalHopfAlgebra \App
\Operad$ defined as follows.  Let $\Reduced$ be the map such that $\Reduced \App w$ is the
subword of $w \in \Operad^*$ consisting of its letters different from $\Unit$. A word $w$ on
$\Operad$ is \Def{reduced} if $\Reduced \App w = w$. Let $\NaturalHopfAlgebra \App \Operad$
be the $\K$-linear span of the set $\Reduced \App \Operad^*$, where $\K$ is any field of
characteristic zero. The bases of $\NaturalHopfAlgebra \App \Operad$ are thus indexed by
$\Reduced \App \Operad^*$, and the \Def{elementary basis} (or \Def{$\BasisE$-basis} for
short) of $\NaturalHopfAlgebra \App \Operad$ is the set $\Bra{\BasisE_w : w \in \Reduced
\App \Operad^*}$. This vector space is endowed with an associative algebra structure through
the product $\Product$ satisfying $\BasisE_w \Product \BasisE_{w'} = \BasisE_{w \Conc w'}$.
Moreover, $\NaturalHopfAlgebra \App \Operad$ is endowed with the coproduct $\Coproduct$
defined as the unique associative algebra morphism satisfying, for any $x \in \Operad$,
\begin{equation} \label{equ:natural_coproduct}
    \Coproduct \App \BasisE_x =
    \sum_{y \in \Operad}
    \;
    \sum_{w \in \Operad^{\Arity \App y}}
    \;
    \Iverson{x = y \Han{w \App 1, \dots, w \App \Length \App w}}
    \;
    \BasisE_{\Reduced \App y} \otimes \BasisE_{\Reduced \App w},
\end{equation}
where the outer $\Iverson{-}$ denotes the Iverson bracket (see \textit{General notations and
conventions}). Due to the fact that $\Operad$ is finitely factorizable,
\eqref{equ:natural_coproduct} is a finite sum. This coproduct endows $\NaturalHopfAlgebra
\App \Operad$ with the structure of a bialgebra. By extending additively $\Degree$ on
$\Operad^*$, the map $\Degree$ defines a grading of $\NaturalHopfAlgebra \App \Operad$.
Thus, $\NaturalHopfAlgebra \App \Operad$ admits an antipode and becomes a Hopf algebra.

For instance, let us consider the \Def{associative operad} $\As$, defined by $\As :=
\Bra{\alpha_n : n \geq 1}$, $\Arity \App \alpha_n := n$, $\alpha_n \Han{\alpha_{m_1}, \dots,
\alpha_{m_n}} := \alpha_{m_1 + \dots + m_n}$, and $\Unit := \alpha_1$. The map $\Degree$
satisfying $\Degree \App \alpha_n = n - 1$ is a grading of $\As$. The bases of
$\NaturalHopfAlgebra \App \As$ are indexed by $\Reduced \App \As^*$ and we have,
\emph{e.g.},
\begin{equation}
    \Coproduct \App \BasisE_{\alpha_4} =
    \BasisE_\epsilon \otimes \BasisE_{\alpha_4}
    + 2 \BasisE_{\alpha_2} \otimes \BasisE_{\alpha_3}
    + \BasisE_{\alpha_2} \otimes \BasisE_{\alpha_2 \alpha_2}
    + 3 \BasisE_{\alpha_3} \otimes \BasisE_{\alpha_2}
    + \BasisE_{\alpha_4} \otimes \BasisE_\epsilon.
\end{equation}
It is shown in~\cite{BG16} that $\NaturalHopfAlgebra \App \As$ is isomorphic to the
\Def{noncommutative Faà di Bruno Hopf algebra} $\NCFdB$ (see~\cite{Foi08}).

\subsection{Polynomial realizations}
For any alphabet $A$, $\K \Angle{A}$ is the $\K$-vector space of \Def{$A$-polynomials},
which are noncommutative polynomials with variables in $A$, having a possibly infinite
support but finite degree. Given two alphabets $A$ and $A'$, let
\begin{math}
    \WordToTensor_{A, A'} : \K \Angle{A \sqcup A'} \to \K \Angle{A} \otimes \K \Angle{A'}
\end{math}
be the linear map such that
\begin{math}
    \WordToTensor_{A, A'} \App w = w_{|A} \otimes w_{|A'}
\end{math}
for any $w \in \Par{A \sqcup A'}^*$.

A \Def{related alphabet} is a set $A$ endowed with a family of $n$-ary relations, $n \geq
1$. As usual, an $n$-ary relation $\NAryRelation$ is a subset of $A^n$, and we denote by
$\NAryRelation\Par{a_1, \dots, a_n}$ the fact that $\Par{a_1, \dots, a_n} \in
\NAryRelation$, and when $n = 2$, we write $a_1 \, \NAryRelation \, a_2$ instead of
$\NAryRelation\Par{a_1, a_2}$. A \Def{polynomial realization} of a Hopf algebra
$\HopfAlgebra$ is a quadruple $\Par{\ClassAlphabets, \AlphabetSum, \Realization_A,
\Alphabet}$ such that
\begin{enumerate}[label=({\sf \roman*})]
    \item $\ClassAlphabets$ is a class of related alphabets;
    \item $\AlphabetSum$ is an associative operation on $\ClassAlphabets$, called
    \Def{disjoint sum};
    \item \label{item:polynomial_realization_3}
    for any related alphabet $A$ of $\ClassAlphabets$, $\Realization_A : \HopfAlgebra \to \K
    \Angle{A}$ is a graded linear map, and for any alphabets $A$ and $A'$ of
    $\ClassAlphabets$ and any $x \in \HopfAlgebra$,
    \begin{math}
        \WordToTensor_{A, A'} \App \Realization_{A \AlphabetSum A'} \App x
        = \Par{\Realization_A \otimes \Realization_{A'}} \circ \Coproduct \App x;
    \end{math}
    \item \label{item:polynomial_realization_4}
    $\Alphabet$ is an alphabet of $\ClassAlphabets$ such that the map
    $\Realization_\Alphabet$ is injective.
\end{enumerate}
Property~\ref{item:polynomial_realization_3} is known as the \Def{alphabet doubling trick}.
This property, enjoyed by polynomial realizations of a large number of Hopf algebras, allows
us to rephrase their coproduct by such alphabet
transformations~\cite{DHT02,NT06,Hiv07,Gir11,FNT14,Foi20}.

\subsection{Terms and free operads}
A \Def{signature} is a graded set $\Signature = \bigsqcup_{n \in \N} \Signature \App n$. A
signature $\Signature$ is \Def{positive} if $\Signature \App 0 = \emptyset$. When all
$\Signature \App n$ are finite, the \Def{profile} of $\Signature$ is the infinite word $w$
such that $w \App i$ is the cardinality of $\Signature \App i - 1$. An
\Def{$\Signature$-term} is either the \Def{leaf} $\Leaf$ or a pair $\Par{\GenG,
\Par{\TermT_1, \dots, \TermT_n}}$ where $\GenG \in \Signature \App n$ and $\TermT_1$, \dots,
$\TermT_n$ are $\Signature$-terms. Thus, an $\Signature$-term is a planar rooted tree where
internal nodes are decorated on $\Signature$. The set of $\Signature$-terms is denoted by
$\SetTerms \App \Signature$. The \Def{free operad} on $\Signature$ is the set $\SetTerms
\App \Signature$ such that the arity $\Arity \App \TermT$ of $\TermT \in \SetTerms \App
\Signature$ is the number of leaves of $\TermT$, for any $\TermT, \TermS_1, \dots,
\TermS_{\Arity \App \TermT} \in \SetTerms \App \Signature$, $\TermT \Han{\TermS_1, \dots,
\TermS_{\Arity \App \TermT}}$ is the $\Signature$-term obtained by grafting the root of each
$\TermS_i$ onto the $i$-th leaf of $\TermT$, and the unit of $\SetTerms \App \Signature$ is
$\Leaf$. The \Def{degree} $\Degree \App \TermT$ of $\TermT \in \SetTerms \App \Signature$ is
the number of internal nodes of $\TermT$. This map $\Degree$ is a grading of $\SetTerms \App
\Signature$.

An \Def{$\Signature$-forest} is any element of $\SetForests \App \Signature :=
\Par{\SetTerms \App \Signature}^*$. We identify each internal node of $\ForestF$ with its
position, starting by $1$, according to its left to right preorder traversal. Let
$\Edge{\ForestF}{j}$ be the binary relation on the set of internal nodes of $\ForestF$ such
that $i \Edge{\ForestF}{j} i'$ holds if $i'$ is the $j$-th child of $i$ in $\ForestF$.
Figure~\ref{fig:example_forest} contains examples of these notions.
\begin{figure}[ht]
    \begin{minipage}{0.42\textwidth}
        \centering
        \begin{equation*}
            \scalebox{.75}{
            \begin{tikzpicture}[Centering,xscale=0.17,yscale=0.4]
                \node[Leaf](0)at(0.00,0.00){};
                \node[Leaf](r)at(0.00,0.75){};
                \draw[Edge](r)--(0);
            \end{tikzpicture}}
            \scalebox{.75}{
            \begin{tikzpicture}[Centering,xscale=0.28,yscale=0.28]
                \node[Leaf](1)at(-1.00,-4.00){};
                \node[Leaf](3)at(1.00,-2.00){};
                \node[Leaf](5)at(2.00,-6.00){};
                \node[Leaf](7)at(4.00,-4.00){};
                \node[Node](0)at(-1.00,-2.00){$\GenA$};
                \node[Node](2)at(1.00,0.00){$\GenC$};
                \node[Node](4)at(2.00,-4.00){$\GenA$};
                \node[Node](6)at(3.00,-2.00){$\GenB$};
                \draw[Edge](0)--(2);
                \draw[Edge](1)--(0);
                \draw[Edge](3)--(2);
                \draw[Edge](4)--(6);
                \draw[Edge](5)--(4);
                \draw[Edge](6)--(2);
                \draw[Edge](7)--(6);
                \node[Leaf](r)at(1.00,2){};
                \draw[Edge](r)--(2);
                \node[left of=2,node distance=10pt,font=\footnotesize,text=ColB]{$1$};
                \node[left of=0,node distance=10pt,font=\footnotesize,text=ColB]{$2$};
                \node[right of=6,node distance=10pt,font=\footnotesize,text=ColB]{$3$};
                \node[left of=4,node distance=10pt,font=\footnotesize,text=ColB]{$4$};
            \end{tikzpicture}}
            \scalebox{.75}{
            \begin{tikzpicture}[Centering,xscale=0.17,yscale=0.4]
                \node[Leaf](0)at(0.00,0.00){};
                \node[Leaf](r)at(0.00,0.75){};
                \draw[Edge](r)--(0);
            \end{tikzpicture}}
            \scalebox{.75}{
            \begin{tikzpicture}[Centering,xscale=0.3,yscale=0.32]
                \node[Leaf](0)at(0.00,-1.50){};
                \node[Leaf](3)at(2.00,-4.50){};
                \node[Leaf](5)at(4.00,-3.00){};
                \node[Node](1)at(1.00,0.00){$\GenB$};
                \node[Node](2)at(2.00,-3.00){$\GenA$};
                \node[Node](4)at(3.00,-1.50){$\GenB$};
                \draw[Edge](0)--(1);
                \draw[Edge](2)--(4);
                \draw[Edge](3)--(2);
                \draw[Edge](4)--(1);
                \draw[Edge](5)--(4);
                \node[Leaf](r)at(1.00,1.5){};
                \draw[Edge](r)--(1);
                \node[left of=1,node distance=10pt,font=\footnotesize,text=ColB]{$5$};
                \node[right of=4,node distance=10pt,font=\footnotesize,text=ColB]{$6$};
                \node[left of=2,node distance=10pt,font=\footnotesize,text=ColB]{$7$};
            \end{tikzpicture}}
            \enspace \xmapsto{\Reduced} \enspace
            \scalebox{.75}{
            \begin{tikzpicture}[Centering,xscale=0.28,yscale=0.25]
                \node[Leaf](1)at(-1.00,-4.00){};
                \node[Leaf](3)at(1.00,-2.00){};
                \node[Leaf](5)at(2.00,-6.00){};
                \node[Leaf](7)at(4.00,-4.00){};
                \node[Node](0)at(-1.00,-2.00){$\GenA$};
                \node[Node](2)at(1.00,0.00){$\GenC$};
                \node[Node](4)at(2.00,-4.00){$\GenA$};
                \node[Node](6)at(3.00,-2.00){$\GenB$};
                \draw[Edge](0)--(2);
                \draw[Edge](1)--(0);
                \draw[Edge](3)--(2);
                \draw[Edge](4)--(6);
                \draw[Edge](5)--(4);
                \draw[Edge](6)--(2);
                \draw[Edge](7)--(6);
                \node[Leaf](r)at(1.00,2){};
                \draw[Edge](r)--(2);
            \end{tikzpicture}}
            \scalebox{.75}{
            \begin{tikzpicture}[Centering,xscale=0.3,yscale=0.29]
                \node[Leaf](0)at(0.00,-1.50){};
                \node[Leaf](3)at(2.00,-4.50){};
                \node[Leaf](5)at(4.00,-3.00){};
                \node[Node](1)at(1.00,0.00){$\GenB$};
                \node[Node](2)at(2.00,-3.00){$\GenA$};
                \node[Node](4)at(3.00,-1.50){$\GenB$};
                \draw[Edge](0)--(1);
                \draw[Edge](2)--(4);
                \draw[Edge](3)--(2);
                \draw[Edge](4)--(1);
                \draw[Edge](5)--(4);
                \node[Leaf](r)at(1.00,1.5){};
                \draw[Edge](r)--(1);
            \end{tikzpicture}}
        \end{equation*}
    \end{minipage}
    \hfill
    \begin{minipage}{0.59\textwidth}
        \caption{A $\SignatureExample$-forest $\ForestF$ and $\Reduced \App \ForestF$, where
        $\SignatureExample$ is the signature $\SignatureExample:= \SignatureExample \App 1
        \sqcup \SignatureExample \App 2 \sqcup \SignatureExample \App 3$ with
        $\SignatureExample \App 1 = \Bra{\GenA}$, $\SignatureExample \App 2 = \Bra{\GenB}$,
        and $\SignatureExample \App 3 = \Bra{\GenC}$. Internal nodes are identified by the
        integers close to them. The degree of $\ForestF$ is $7$, and we have for instance $1
        \Edge{\ForestF}{1} 2$, $1 \Edge{\ForestF}{3} 3$, and $5 \Edge{\ForestF}{2} 6$.}
        \label{fig:example_forest}
    \end{minipage}
\end{figure}

Note that $\NaturalHopfAlgebra \App \SetTerms \App \Signature$ is a Hopf algebra graded by
$\Degree$ and whose bases are indexed by the set $\Reduced \App \SetForests \App \Signature$
of reduced $\Signature$-forests.

\section{Lattices on forests and alternative bases} \label{sec:alternative_bases}
Let $\Signature$ be a positive signature and $\EasterlyWindCovering$ be the binary relation
on $\Reduced \App \SetForests \App \Signature$ such that $\ForestF \EasterlyWindCovering
\ForestF'$ if~$\ForestF'$ can be obtained from $\ForestF$ by selecting an internal node $i$,
detaching the subterm rooted at~$i$, and grafting it to the rightmost leaf of its former
immediate left brother. For instance,
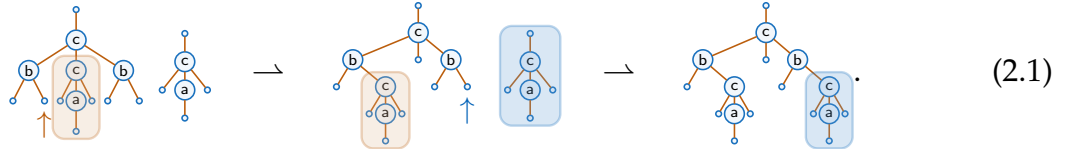
\begin{equation}
    \scalebox{.75}{
    \begin{tikzpicture}[Centering,xscale=0.28,yscale=0.18]
        \node[Leaf](0)at(0.00,-6.00){};
        \node[Leaf](11)at(8.00,-6.00){};
        \node[Leaf](2)at(2.00,-6.00){};
        \node[Leaf](4)at(3.00,-6.00){};
        \node[Leaf](7)at(4.00,-9.00){};
        \node[Leaf](8)at(5.00,-6.00){};
        \node[Leaf](9)at(6.00,-6.00){};
        \node[Node](1)at(1.00,-3.00){$\GenB$};
        \node[Node](10)at(7.00,-3.00){$\GenB$};
        \node[Node](3)at(4.00,0.00){$\GenC$};
        \node[Node](5)at(4.00,-3.00){$\GenC$};
        \node[Node](6)at(4.00,-6.00){$\GenA$};
        \draw[Edge](0)--(1);
        \draw[Edge](1)--(3);
        \draw[Edge](10)--(3);
        \draw[Edge](11)--(10);
        \draw[Edge](2)--(1);
        \draw[Edge](4)--(5);
        \draw[Edge](5)--(3);
        \draw[Edge](6)--(5);
        \draw[Edge](7)--(6);
        \draw[Edge](8)--(5);
        \draw[Edge](9)--(10);
        \node[Leaf](r)at(4.00,3){};
        \draw[Edge](r)--(3);
        \node[FitNodes,fit=(4)(5)(6)(7)(8)]{};
        \node[below of=2,node distance=12pt,font=\large,text=ColB]{$\uparrow$};
    \end{tikzpicture}}
    \;
    \scalebox{.75}{
    \begin{tikzpicture}[Centering,xscale=0.4,yscale=0.3]
        \node[Leaf](0)at(0.00,-1.67){};
        \node[Leaf](3)at(1.00,-3.33){};
        \node[Leaf](4)at(2.00,-1.67){};
        \node[Node](1)at(1.00,0.00){$\GenC$};
        \node[Node](2)at(1.00,-1.67){$\GenA$};
        \draw[Edge](0)--(1);
        \draw[Edge](2)--(1);
        \draw[Edge](3)--(2);
        \draw[Edge](4)--(1);
        \node[Leaf](r)at(1.00,1.67){};
        \draw[Edge](r)--(1);
    \end{tikzpicture}}
    \quad \EasterlyWindCovering \quad
    \scalebox{.75}{
    \begin{tikzpicture}[Centering,xscale=0.29,yscale=0.2]
        \node[Leaf](0)at(0.00,-4.80){};
        \node[Leaf](11)at(8.00,-4.80){};
        \node[Leaf](2)at(2.00,-7.20){};
        \node[Leaf](5)at(3.00,-9.60){};
        \node[Leaf](6)at(4.00,-7.20){};
        \node[Leaf](8)at(5.00,-2.40){};
        \node[Leaf](9)at(6.00,-4.80){};
        \node[Node](1)at(1.00,-2.40){$\GenB$};
        \node[Node](10)at(7.00,-2.40){$\GenB$};
        \node[Node](3)at(3.00,-4.80){$\GenC$};
        \node[Node](4)at(3.00,-7.20){$\GenA$};
        \node[Node](7)at(5.00,0.00){$\GenC$};
        \draw[Edge](0)--(1);
        \draw[Edge](1)--(7);
        \draw[Edge](10)--(7);
        \draw[Edge](11)--(10);
        \draw[Edge](2)--(3);
        \draw[Edge](3)--(1);
        \draw[Edge](4)--(3);
        \draw[Edge](5)--(4);
        \draw[Edge](6)--(3);
        \draw[Edge](8)--(7);
        \draw[Edge](9)--(10);
        \node[Leaf](r)at(5.00,2.40){};
        \draw[Edge](r)--(7);
        \node[FitNodes,fit=(2)(3)(4)(5)(6)]{};
        \node[below of=11,node distance=12pt,font=\large,text=ColA]{$\uparrow$};
    \end{tikzpicture}}
    \;
    \scalebox{.75}{
    \begin{tikzpicture}[Centering,xscale=0.4,yscale=0.3]
        \node[Leaf](0)at(0.00,-1.67){};
        \node[Leaf](3)at(1.00,-3.33){};
        \node[Leaf](4)at(2.00,-1.67){};
        \node[Node](1)at(1.00,0.00){$\GenC$};
        \node[Node](2)at(1.00,-1.67){$\GenA$};
        \draw[Edge](0)--(1);
        \draw[Edge](2)--(1);
        \draw[Edge](3)--(2);
        \draw[Edge](4)--(1);
        \node[Leaf](r)at(1.00,1.67){};
        \draw[Edge](r)--(1);
        \node[FitNodes,draw=ColA,fill=ColA!60,fit=(0)(1)(2)(3)(4)(r)]{};
    \end{tikzpicture}}
    \quad \EasterlyWindCovering \quad
    \scalebox{.75}{
    \begin{tikzpicture}[Centering,xscale=0.28,yscale=0.15]
        \node[Leaf](0)at(0.00,-6.40){};
        \node[Leaf](11)at(8.00,-9.60){};
        \node[Leaf](14)at(9.00,-12.80){};
        \node[Leaf](15)at(10.00,-9.60){};
        \node[Leaf](2)at(2.00,-9.60){};
        \node[Leaf](5)at(3.00,-12.80){};
        \node[Leaf](6)at(4.00,-9.60){};
        \node[Leaf](8)at(5.00,-3.20){};
        \node[Leaf](9)at(6.00,-6.40){};
        \node[Node](1)at(1.00,-3.20){$\GenB$};
        \node[Node](10)at(7.00,-3.20){$\GenB$};
        \node[Node](12)at(9.00,-6.40){$\GenC$};
        \node[Node](13)at(9.00,-9.60){$\GenA$};
        \node[Node](3)at(3.00,-6.40){$\GenC$};
        \node[Node](4)at(3.00,-9.60){$\GenA$};
        \node[Node](7)at(5.00,0.00){$\GenC$};
        \draw[Edge](0)--(1);
        \draw[Edge](1)--(7);
        \draw[Edge](10)--(7);
        \draw[Edge](11)--(12);
        \draw[Edge](12)--(10);
        \draw[Edge](13)--(12);
        \draw[Edge](14)--(13);
        \draw[Edge](15)--(12);
        \draw[Edge](2)--(3);
        \draw[Edge](3)--(1);
        \draw[Edge](4)--(3);
        \draw[Edge](5)--(4);
        \draw[Edge](6)--(3);
        \draw[Edge](8)--(7);
        \draw[Edge](9)--(10);
        \node[Leaf](r)at(5.00,3.20){};
        \draw[Edge](r)--(7);
        \node[FitNodes,draw=ColA,fill=ColA!60,fit=(11)(12)(13)(14)(15)]{};
    \end{tikzpicture}}.
\end{equation}
Let $\Leq$ be the reflexive and transitive closure of $\EasterlyWindCovering$.

\begin{Proposition} \label{prop:easterly_wind_poset}
    For any positive signature $\Signature$, the relation $\Leq$ is a partial order on
    $\Reduced \App \SetForests \App \Signature$. Moreover, $\EasterlyWindCovering$ is the
    covering relation of this partial order.
\end{Proposition}

The poset $\Par{\Reduced \App \SetForests \App \Signature, \Leq}$ is called the
\Def{$\Signature$-easterly wind poset}. This name is derived from the observation that the
covering relation of this poset involves detaching a subterm from the east and then
attaching it to the west, as if an easterly breeze is blowing on the forest. Let $\Corolla :
\Signature \to \SetTerms \App \Signature$ be the map such that $\Corolla \App \GenG$ is the
$\Signature$-term of degree $1$ whose root is decorated by $\GenG$. Given $w \in
\Signature^*$, let $\EasterlyWindPoset \App w$ be the interval of the $\Signature$-easterly
wind poset whose lower bound is the $\Signature$-forest $\ForestF_{\min}$ of length $\Length
\App w$ such that $\ForestF_{\min} \App i := \Corolla \App w \App i$ for any $i \in [\Length
\App w]$, and the upper bound is the $\Signature$-forest $\ForestF_{\max}$ of length $1$
such that
\begin{math}
    \ForestF_{\max} \App 1
    := (\Corolla \App w \App 1) \circ_1 \dots \circ_1 (\Corolla \App w \App \Length \App w),
\end{math}
where $\circ_i$ is the partial composition map of the free operad~$\SetTerms \App
\Signature$. Figure~\ref{fig:example_easterly_wind_poset} shows an example of an
$\Signature$-easterly wind lattice.
\begin{figure}[ht]
    \begin{minipage}{0.6\textwidth}
        \centering
        \scalebox{0.67}{
        \begin{tikzpicture}[Centering,xscale=2.5,yscale=1.88]
            \node[GraphLabeledVertex](c000a0b00)at(0,0){
                \begin{tikzpicture}[Centering,xscale=0.2,yscale=0.24]
                    \node[Leaf](0)at(0.00,-2.00){};
                    \node[Leaf](2)at(1.00,-2.00){};
                    \node[Leaf](3)at(2.00,-2.00){};
                    \node[Node](1)at(1.00,0.00){$\GenC$};
                    \draw[Edge](0)--(1);
                    \draw[Edge](2)--(1);
                    \draw[Edge](3)--(1);
                    \node[Leaf](r)at(1.00,2){};
                    \draw[Edge](r)--(1);
                \end{tikzpicture}
                \begin{tikzpicture}[Centering,xscale=0.17,yscale=0.4]
                    \node[Leaf](1)at(0.00,-1.00){};
                    \node[Node](0)at(0.00,0.00){$\GenA$};
                    \draw[Edge](1)--(0);
                    \node[Leaf](r)at(0.00,1){};
                    \draw[Edge](r)--(0);
                \end{tikzpicture}
                \begin{tikzpicture}[Centering,xscale=0.2,yscale=0.3]
                    \node[Leaf](0)at(0.00,-1.50){};
                    \node[Leaf](2)at(2.00,-1.50){};
                    \node[Node](1)at(1.00,0.00){$\GenB$};
                    \draw[Edge](0)--(1);
                    \draw[Edge](2)--(1);
                    \node[Leaf](r)at(1.00,1.5){};
                    \draw[Edge](r)--(1);
                \end{tikzpicture}
            };
            \node[GraphLabeledVertex](c00a0b00)at(1,-1){
                \begin{tikzpicture}[Centering,xscale=0.20,yscale=0.28]
                    \node[Leaf](0)at(-1.00,-1.67){};
                    \node[Leaf](2)at(1.00,-1.67){};
                    \node[Leaf](4)at(3.00,-3.33){};
                    \node[Node](1)at(1.00,0.00){$\GenC$};
                    \node[Node](3)at(3.00,-1.67){$\GenA$};
                    \draw[Edge](0)--(1);
                    \draw[Edge](2)--(1);
                    \draw[Edge](3)--(1);
                    \draw[Edge](4)--(3);
                    \node[Leaf](r)at(1.00,1.67){};
                    \draw[Edge](r)--(1);
                \end{tikzpicture}
                \begin{tikzpicture}[Centering,xscale=0.2,yscale=0.3]
                    \node[Leaf](0)at(0.00,-1.50){};
                    \node[Leaf](2)at(2.00,-1.50){};
                    \node[Node](1)at(1.00,0.00){$\GenB$};
                    \draw[Edge](0)--(1);
                    \draw[Edge](2)--(1);
                    \node[Leaf](r)at(1.00,1.5){};
                    \draw[Edge](r)--(1);
                \end{tikzpicture}
            };
            \node[GraphLabeledVertex](c00ab00)at(-1,-3){
                \begin{tikzpicture}[Centering,xscale=0.24,yscale=0.28]
                    \node[Leaf](0)at(0.00,-1.75){};
                    \node[Leaf](2)at(1.00,-1.75){};
                    \node[Leaf](4)at(2.00,-5.25){};
                    \node[Leaf](6)at(4.00,-5.25){};
                    \node[Node](1)at(1.00,0.00){$\GenC$};
                    \node[Node](3)at(3.00,-1.75){$\GenA$};
                    \node[Node](5)at(3.00,-3.50){$\GenB$};
                    \draw[Edge](0)--(1);
                    \draw[Edge](2)--(1);
                    \draw[Edge](3)--(1);
                    \draw[Edge](4)--(5);
                    \draw[Edge](5)--(3);
                    \draw[Edge](6)--(5);
                    \node[Leaf](r)at(1.00,1.75){};
                    \draw[Edge](r)--(1);
                \end{tikzpicture}
            };
            \node[GraphLabeledVertex](c000ab00)at(-2,-2){
                \begin{tikzpicture}[Centering,xscale=0.2,yscale=0.24]
                    \node[Leaf](0)at(0.00,-2.00){};
                    \node[Leaf](2)at(1.00,-2.00){};
                    \node[Leaf](3)at(2.00,-2.00){};
                    \node[Node](1)at(1.00,0.00){$\GenC$};
                    \draw[Edge](0)--(1);
                    \draw[Edge](2)--(1);
                    \draw[Edge](3)--(1);
                    \node[Leaf](r)at(1.00,2){};
                    \draw[Edge](r)--(1);
                \end{tikzpicture}
                \begin{tikzpicture}[Centering,xscale=0.2,yscale=0.38]
                    \node[Leaf](1)at(0.00,-2.67){};
                    \node[Leaf](3)at(2.00,-2.67){};
                    \node[Node](0)at(1.00,0.00){$\GenA$};
                    \node[Node](2)at(1.00,-1.33){$\GenB$};
                    \draw[Edge](1)--(2);
                    \draw[Edge](2)--(0);
                    \draw[Edge](3)--(2);
                    \node[Leaf](r)at(1.00,1.33){};
                    \draw[Edge](r)--(0);
                \end{tikzpicture}
            };
            \node[GraphLabeledVertex](c0a00b00)at(2,-2){
                \begin{tikzpicture}[Centering,xscale=0.24,yscale=0.28]
                    \node[Leaf](0)at(-1.00,-1.67){};
                    \node[Leaf](3)at(1.00,-3.33){};
                    \node[Leaf](4)at(3.00,-1.67){};
                    \node[Node](1)at(1.00,0.00){$\GenC$};
                    \node[Node](2)at(1.00,-1.67){$\GenA$};
                    \draw[Edge](0)--(1);
                    \draw[Edge](2)--(1);
                    \draw[Edge](3)--(2);
                    \draw[Edge](4)--(1);
                    \node[Leaf](r)at(1.00,1.67){};
                    \draw[Edge](r)--(1);
                \end{tikzpicture}
                \begin{tikzpicture}[Centering,xscale=0.2,yscale=0.3]
                    \node[Leaf](0)at(0.00,-1.50){};
                    \node[Leaf](2)at(2.00,-1.50){};
                    \node[Node](1)at(1.00,0.00){$\GenB$};
                    \draw[Edge](0)--(1);
                    \draw[Edge](2)--(1);
                    \node[Leaf](r)at(1.00,1.5){};
                    \draw[Edge](r)--(1);
                \end{tikzpicture}
            };
            \node[GraphLabeledVertex](ca000b00)at(3,-3){
                \begin{tikzpicture}[Centering,xscale=0.20,yscale=0.28]
                    \node[Leaf](1)at(-1.00,-3.33){};
                    \node[Leaf](3)at(1.00,-1.67){};
                    \node[Leaf](4)at(3.00,-1.67){};
                    \node[Node](0)at(-1.00,-1.67){$\GenA$};
                    \node[Node](2)at(1.00,0.00){$\GenC$};
                    \draw[Edge](0)--(2);
                    \draw[Edge](1)--(0);
                    \draw[Edge](3)--(2);
                    \draw[Edge](4)--(2);
                    \node[Leaf](r)at(1.00,1.67){};
                    \draw[Edge](r)--(2);
                \end{tikzpicture}
                \begin{tikzpicture}[Centering,xscale=0.2,yscale=0.3]
                    \node[Leaf](0)at(0.00,-1.50){};
                    \node[Leaf](2)at(2.00,-1.50){};
                    \node[Node](1)at(1.00,0.00){$\GenB$};
                    \draw[Edge](0)--(1);
                    \draw[Edge](2)--(1);
                    \node[Leaf](r)at(1.00,1.5){};
                    \draw[Edge](r)--(1);
                \end{tikzpicture}
            };
            \node[GraphLabeledVertex](c0a0b00)at(1,-3){
                \begin{tikzpicture}[Centering,xscale=0.24,yscale=0.21]
                    \node[Leaf](0)at(-1.00,-2.33){};
                    \node[Leaf](3)at(1.00,-4.67){};
                    \node[Leaf](4)at(2.00,-4.67){};
                    \node[Leaf](6)at(4.00,-4.67){};
                    \node[Node](1)at(1.00,0.00){$\GenC$};
                    \node[Node](2)at(1.00,-2.33){$\GenA$};
                    \node[Node](5)at(3.00,-2.33){$\GenB$};
                    \draw[Edge](0)--(1);
                    \draw[Edge](2)--(1);
                    \draw[Edge](3)--(2);
                    \draw[Edge](4)--(5);
                    \draw[Edge](5)--(1);
                    \draw[Edge](6)--(5);
                    \node[Leaf](r)at(1.00,2.33){};
                    \draw[Edge](r)--(1);
                \end{tikzpicture}
            };
            \node[GraphLabeledVertex](ca00b00)at(2,-4){
                \begin{tikzpicture}[Centering,xscale=0.24,yscale=0.21]
                    \node[Leaf](1)at(-1.00,-4.67){};
                    \node[Leaf](3)at(1.00,-2.33){};
                    \node[Leaf](4)at(2.00,-4.67){};
                    \node[Leaf](6)at(4.00,-4.67){};
                    \node[Node](0)at(-1.00,-2.33){$\GenA$};
                    \node[Node](2)at(1.00,0.00){$\GenC$};
                    \node[Node](5)at(3.00,-2.33){$\GenB$};
                    \draw[Edge](0)--(2);
                    \draw[Edge](1)--(0);
                    \draw[Edge](3)--(2);
                    \draw[Edge](4)--(5);
                    \draw[Edge](5)--(2);
                    \draw[Edge](6)--(5);
                    \node[Leaf](r)at(1.00,2.33){};
                    \draw[Edge](r)--(2);
                \end{tikzpicture}
            };
            \node[GraphLabeledVertex](c0ab000)at(0,-4){
                \begin{tikzpicture}[Centering,xscale=0.24,yscale=0.28]
                    \node[Leaf](0)at(0.00,-1.75){};
                    \node[Leaf](3)at(1.00,-5.25){};
                    \node[Leaf](5)at(3.00,-5.25){};
                    \node[Leaf](6)at(4.00,-1.75){};
                    \node[Node](1)at(2.00,0.00){$\GenC$};
                    \node[Node](2)at(2.00,-1.75){$\GenA$};
                    \node[Node](4)at(2.00,-3.50){$\GenB$};
                    \draw[Edge](0)--(1);
                    \draw[Edge](2)--(1);
                    \draw[Edge](3)--(4);
                    \draw[Edge](4)--(2);
                    \draw[Edge](5)--(4);
                    \draw[Edge](6)--(1);
                    \node[Leaf](r)at(2.00,1.75){};
                    \draw[Edge](r)--(1);
                \end{tikzpicture}
            };
            \node[GraphLabeledVertex](ca0b000)at(3,-5){
                \begin{tikzpicture}[Centering,xscale=0.24,yscale=0.20]
                    \node[Leaf](1)at(0.00,-4.67){};
                    \node[Leaf](3)at(1.00,-4.67){};
                    \node[Leaf](5)at(3.00,-4.67){};
                    \node[Leaf](6)at(4.00,-2.33){};
                    \node[Node](0)at(0.00,-2.33){$\GenA$};
                    \node[Node](2)at(2.00,0.00){$\GenC$};
                    \node[Node](4)at(2.00,-2.33){$\GenB$};
                    \draw[Edge](0)--(2);
                    \draw[Edge](1)--(0);
                    \draw[Edge](3)--(4);
                    \draw[Edge](4)--(2);
                    \draw[Edge](5)--(4);
                    \draw[Edge](6)--(2);
                    \node[Leaf](r)at(2.00,2.33){};
                    \draw[Edge](r)--(2);
                \end{tikzpicture}
            };
            \node[GraphLabeledVertex](cab0000)at(2,-6){
                \begin{tikzpicture}[Centering,xscale=0.24,yscale=0.28]
                    \node[Leaf](1)at(0.00,-5.25){};
                    \node[Leaf](3)at(2.00,-5.25){};
                    \node[Leaf](5)at(3.00,-1.75){};
                    \node[Leaf](6)at(4.00,-1.75){};
                    \node[Node](0)at(1.00,-1.75){$\GenA$};
                    \node[Node](2)at(1.00,-3.50){$\GenB$};
                    \node[Node](4)at(3.00,0.00){$\GenC$};
                    \draw[Edge](0)--(4);
                    \draw[Edge](1)--(2);
                    \draw[Edge](2)--(0);
                    \draw[Edge](3)--(2);
                    \draw[Edge](5)--(4);
                    \draw[Edge](6)--(4);
                    \node[Leaf](r)at(3.00,1.75){};
                    \draw[Edge](r)--(4);
                \end{tikzpicture}
            };
            \draw[GraphEdge](c000a0b00)--(c00a0b00);
            \draw[GraphEdge](c00a0b00)--(c00ab00);
            \draw[GraphEdge](c000a0b00)--(c000ab00);
            \draw[GraphEdge](c000ab00)--(c00ab00);
            \draw[GraphEdge](c00a0b00)--(c0a00b00);
            \draw[GraphEdge](c0a00b00)--(ca000b00);
            \draw[GraphEdge](ca000b00)--(ca00b00);
            \draw[GraphEdge](c0a00b00)--(c0a0b00);
            \draw[GraphEdge](c00ab00)--(c0ab000);
            \draw[GraphEdge](c0a0b00)--(c0ab000);
            \draw[GraphEdge](c0a0b00)--(ca00b00);
            \draw[GraphEdge](ca00b00)--(ca0b000);
            \draw[GraphEdge](ca0b000)--(cab0000);
            \draw[GraphEdge](c0ab000)--(cab0000);
        \end{tikzpicture}}
    \end{minipage}
    \hfill
    \begin{minipage}{0.4\textwidth}
        \caption{The Hasse diagram of the poset $\EasterlyWindPoset \App \GenC \GenA \GenB$
        where $\GenA$, $\GenB$, and $\GenC$ belong to the signature $\SignatureExample$
        defined in Figure~\ref{fig:example_forest}. \\[1ex]
        The reduced $\SignatureExample$-forest at the top is the minimal element of the
        poset, and the one at the bottom is the maximal element.}
        \label{fig:example_easterly_wind_poset}
    \end{minipage}
\end{figure}

\begin{Theorem} \label{thm:easterly_wind_poset_lattice}
    For any positive signature $\Signature$ and any $w \in \Signature^*$, the poset
    $\EasterlyWindPoset \App w$ is a lattice.
\end{Theorem}

Let us now use these lattices to build new bases of $\NaturalHopfAlgebra \App \SetTerms
\App \Signature$. For any $\ForestF \in \Reduced \App \SetForests \App \Signature$, let
\begin{equation}
    \BasisF_\ForestF :=
    \sum_{\ForestF' \in \Reduced \App \SetForests \App \Signature}
    \Iverson{\ForestF \Leq \ForestF'} \;
    \Mobius_{\Leq}\Par{\ForestF, \ForestF'} \;
    \BasisE_{\ForestF'}
    \quad \text{and} \quad
    \BasisH_\ForestF :=
    \sum_{\ForestF' \in \Reduced \App \SetForests \App \Signature}
    \Iverson{\ForestF' \Leq \ForestF} \;
    \BasisF_{\ForestF'},
\end{equation}
where $\Mobius_{\Leq}$ is the Möbius function of the $\Signature$-easterly wind poset. For
instance, in $\NaturalHopfAlgebra \App \SetTerms \App \SignatureExample$,
\begin{equation}
    \BasisF_{
        \scalebox{.75}{
        \begin{tikzpicture}[Centering,xscale=0.34,yscale=0.28]
            \node[Leaf](0)at(0.00,-1.67){};
            \node[Leaf](3)at(1.00,-3.33){};
            \node[Leaf](4)at(2.00,-1.67){};
            \node[Node](1)at(1.00,0.00){$\GenC$};
            \node[Node](2)at(1.00,-1.67){$\GenA$};
            \draw[Edge](0)--(1);
            \draw[Edge](2)--(1);
            \draw[Edge](3)--(2);
            \draw[Edge](4)--(1);
            \node[Leaf](r)at(1.00,1.67){};
            \draw[Edge](r)--(1);
        \end{tikzpicture}}
        \scalebox{.75}{
        \begin{tikzpicture}[Centering,xscale=0.2,yscale=0.3]
            \node[Leaf](0)at(0.00,-1.50){};
            \node[Leaf](2)at(2.00,-1.50){};
            \node[Node](1)at(1.00,0.00){$\GenB$};
            \draw[Edge](0)--(1);
            \draw[Edge](2)--(1);
            \node[Leaf](r)at(1.00,1.5){};
            \draw[Edge](r)--(1);
        \end{tikzpicture}}}
    =
    \BasisE_{
        \scalebox{.75}{
        \begin{tikzpicture}[Centering,xscale=0.34,yscale=0.28]
            \node[Leaf](0)at(0.00,-1.67){};
            \node[Leaf](3)at(1.00,-3.33){};
            \node[Leaf](4)at(2.00,-1.67){};
            \node[Node](1)at(1.00,0.00){$\GenC$};
            \node[Node](2)at(1.00,-1.67){$\GenA$};
            \draw[Edge](0)--(1);
            \draw[Edge](2)--(1);
            \draw[Edge](3)--(2);
            \draw[Edge](4)--(1);
            \node[Leaf](r)at(1.00,1.5){};
            \draw[Edge](r)--(1);
        \end{tikzpicture}}
        \scalebox{.75}{
        \begin{tikzpicture}[Centering,xscale=0.2,yscale=0.3]
            \node[Leaf](0)at(0.00,-1.50){};
            \node[Leaf](2)at(2.00,-1.50){};
            \node[Node](1)at(1.00,0.00){$\GenB$};
            \draw[Edge](0)--(1);
            \draw[Edge](2)--(1);
            \node[Leaf](r)at(1.00,1.5){};
            \draw[Edge](r)--(1);
        \end{tikzpicture}}}
    -
    \BasisE_{
        \scalebox{.75}{
        \begin{tikzpicture}[Centering,xscale=0.24,yscale=0.21]
            \node[Leaf](0)at(-1.00,-2.33){};
            \node[Leaf](3)at(1.00,-4.67){};
            \node[Leaf](4)at(2.00,-4.67){};
            \node[Leaf](6)at(4.00,-4.67){};
            \node[Node](1)at(1.00,0.00){$\GenC$};
            \node[Node](2)at(1.00,-2.33){$\GenA$};
            \node[Node](5)at(3.00,-2.33){$\GenB$};
            \draw[Edge](0)--(1);
            \draw[Edge](2)--(1);
            \draw[Edge](3)--(2);
            \draw[Edge](4)--(5);
            \draw[Edge](5)--(1);
            \draw[Edge](6)--(5);
            \node[Leaf](r)at(1.00,2.33){};
            \draw[Edge](r)--(1);
        \end{tikzpicture}}}
    -
    \BasisE_{
        \scalebox{.75}{
        \begin{tikzpicture}[Centering,xscale=0.22,yscale=0.28]
            \node[Leaf](1)at(-0.50,-3.33){};
            \node[Leaf](3)at(1.00,-1.67){};
            \node[Leaf](4)at(2.50,-1.67){};
            \node[Node](0)at(-0.50,-1.67){$\GenA$};
            \node[Node](2)at(1.00,0.00){$\GenC$};
            \draw[Edge](0)--(2);
            \draw[Edge](1)--(0);
            \draw[Edge](3)--(2);
            \draw[Edge](4)--(2);
            \node[Leaf](r)at(1.00,1.67){};
            \draw[Edge](r)--(2);
        \end{tikzpicture}}
        \scalebox{.75}{
        \begin{tikzpicture}[Centering,xscale=0.2,yscale=0.3]
            \node[Leaf](0)at(0.00,-1.50){};
            \node[Leaf](2)at(2.00,-1.50){};
            \node[Node](1)at(1.00,0.00){$\GenB$};
            \draw[Edge](0)--(1);
            \draw[Edge](2)--(1);
            \node[Leaf](r)at(1.00,1.5){};
            \draw[Edge](r)--(1);
        \end{tikzpicture}}}
    +
    \BasisE_{
        \scalebox{.75}{
        \begin{tikzpicture}[Centering,xscale=0.20,yscale=0.21]
            \node[Leaf](1)at(-1.00,-4.67){};
            \node[Leaf](3)at(1.00,-2.33){};
            \node[Leaf](4)at(2.00,-4.67){};
            \node[Leaf](6)at(4.00,-4.67){};
            \node[Node](0)at(-1.00,-2.33){$\GenA$};
            \node[Node](2)at(1.00,0.00){$\GenC$};
            \node[Node](5)at(3.00,-2.33){$\GenB$};
            \draw[Edge](0)--(2);
            \draw[Edge](1)--(0);
            \draw[Edge](3)--(2);
            \draw[Edge](4)--(5);
            \draw[Edge](5)--(2);
            \draw[Edge](6)--(5);
            \node[Leaf](r)at(1.00,2.33){};
            \draw[Edge](r)--(2);
        \end{tikzpicture}}},
\end{equation}
\begin{equation}
    \BasisH_{
        \scalebox{.75}{
        \begin{tikzpicture}[Centering,xscale=0.34,yscale=0.28]
            \node[Leaf](0)at(0.00,-1.67){};
            \node[Leaf](3)at(1.00,-3.33){};
            \node[Leaf](4)at(2.00,-1.67){};
            \node[Node](1)at(1.00,0.00){$\GenC$};
            \node[Node](2)at(1.00,-1.67){$\GenA$};
            \draw[Edge](0)--(1);
            \draw[Edge](2)--(1);
            \draw[Edge](3)--(2);
            \draw[Edge](4)--(1);
            \node[Leaf](r)at(1.00,1.67){};
            \draw[Edge](r)--(1);
        \end{tikzpicture}}
        \scalebox{.75}{
        \begin{tikzpicture}[Centering,xscale=0.2,yscale=0.24]
            \node[Leaf](0)at(0.00,-1.50){};
            \node[Leaf](2)at(2.00,-1.50){};
            \node[Node](1)at(1.00,0.00){$\GenB$};
            \draw[Edge](0)--(1);
            \draw[Edge](2)--(1);
            \node[Leaf](r)at(1.00,1.5){};
            \draw[Edge](r)--(1);
        \end{tikzpicture}}}
    =
    \BasisF_{
        \scalebox{.75}{
        \begin{tikzpicture}[Centering,xscale=0.2,yscale=0.22]
            \node[Leaf](0)at(0.00,-2.00){};
            \node[Leaf](2)at(1.00,-2.00){};
            \node[Leaf](3)at(2.00,-2.00){};
            \node[Node](1)at(1.00,0.00){$\GenC$};
            \draw[Edge](0)--(1);
            \draw[Edge](2)--(1);
            \draw[Edge](3)--(1);
            \node[Leaf](r)at(1.00,2){};
            \draw[Edge](r)--(1);
        \end{tikzpicture}}
        \scalebox{.75}{
        \begin{tikzpicture}[Centering,xscale=0.17,yscale=0.38]
            \node[Leaf](1)at(0.00,-1.00){};
            \node[Node](0)at(0.00,0.00){$\GenA$};
            \draw[Edge](1)--(0);
            \node[Leaf](r)at(0.00,1){};
            \draw[Edge](r)--(0);
        \end{tikzpicture}}
        \scalebox{.75}{
        \begin{tikzpicture}[Centering,xscale=0.2,yscale=0.24]
            \node[Leaf](0)at(0.00,-1.50){};
            \node[Leaf](2)at(2.00,-1.50){};
            \node[Node](1)at(1.00,0.00){$\GenB$};
            \draw[Edge](0)--(1);
            \draw[Edge](2)--(1);
            \node[Leaf](r)at(1.00,1.5){};
            \draw[Edge](r)--(1);
        \end{tikzpicture}}}
    +
    \BasisF_{
        \scalebox{.75}{
        \begin{tikzpicture}[Centering,xscale=0.22,yscale=0.28]
            \node[Leaf](0)at(-0.5,-1.67){};
            \node[Leaf](2)at(1.00,-1.67){};
            \node[Leaf](4)at(2.5,-3.33){};
            \node[Node](1)at(1.00,0.00){$\GenC$};
            \node[Node](3)at(2.5,-1.67){$\GenA$};
            \draw[Edge](0)--(1);
            \draw[Edge](2)--(1);
            \draw[Edge](3)--(1);
            \draw[Edge](4)--(3);
            \node[Leaf](r)at(1.00,1.67){};
            \draw[Edge](r)--(1);
        \end{tikzpicture}}
        \scalebox{.75}{
        \begin{tikzpicture}[Centering,xscale=0.2,yscale=0.24]
            \node[Leaf](0)at(0.00,-1.50){};
            \node[Leaf](2)at(2.00,-1.50){};
            \node[Node](1)at(1.00,0.00){$\GenB$};
            \draw[Edge](0)--(1);
            \draw[Edge](2)--(1);
            \node[Leaf](r)at(1.00,1.5){};
            \draw[Edge](r)--(1);
        \end{tikzpicture}}}
    +
    \BasisF_{
        \scalebox{.75}{
        \begin{tikzpicture}[Centering,xscale=0.34,yscale=0.28]
            \node[Leaf](0)at(0.00,-1.67){};
            \node[Leaf](3)at(1.00,-3.33){};
            \node[Leaf](4)at(2.00,-1.67){};
            \node[Node](1)at(1.00,0.00){$\GenC$};
            \node[Node](2)at(1.00,-1.67){$\GenA$};
            \draw[Edge](0)--(1);
            \draw[Edge](2)--(1);
            \draw[Edge](3)--(2);
            \draw[Edge](4)--(1);
            \node[Leaf](r)at(1.00,1.67){};
            \draw[Edge](r)--(1);
        \end{tikzpicture}}
        \scalebox{.75}{
        \begin{tikzpicture}[Centering,xscale=0.2,yscale=0.24]
            \node[Leaf](0)at(0.00,-1.50){};
            \node[Leaf](2)at(2.00,-1.50){};
            \node[Node](1)at(1.00,0.00){$\GenB$};
            \draw[Edge](0)--(1);
            \draw[Edge](2)--(1);
            \node[Leaf](r)at(1.00,1.5){};
            \draw[Edge](r)--(1);
        \end{tikzpicture}}}.
\end{equation}
By Möbius inversion and triangularity,
\begin{equation}
    \BasisE_\ForestF =
    \sum_{\ForestF' \in \Reduced \App \SetForests \App \Signature}
    \Iverson{\ForestF \Leq \ForestF'} \; \BasisF_{\ForestF'}
    \quad \text{and} \quad
    \BasisF_\ForestF =
    \sum_{\ForestF' \in \Reduced \App \SetForests \App \Signature}
    \Iverson{\ForestF' \Leq \ForestF} \;
    \Mobius_{\Leq}\Par{\ForestF', \ForestF} \;
    \BasisH_{\ForestF'},
\end{equation}
so that the sets $\Bra{\BasisF_\ForestF : \ForestF \in \Reduced \App \SetForests \App
\Signature}$ and $\Bra{\BasisH_\ForestF : \ForestF \in \Reduced \App \SetForests \App
\Signature}$ are bases of $\NaturalHopfAlgebra \App \SetTerms \App \Signature$, called
respectively the \Def{fundamental basis} (or \Def{$\BasisF$-basis} for short) and the
\Def{homogeneous basis} (or \Def{$\BasisH$-basis} for short).

The $i$-th leaf of $\TermT \in \SetTerms \App \Signature$ is \Def{extremal} if, by setting
$\Signature' := \Signature \sqcup \{\DummyConstant\}$ where $\DummyConstant$ is a constant
of arity $0$, the only internal node of the $\Signature'$-term $\TermT \circ_i
\DummyConstant$ which is decorated by $\DummyConstant$ is the last one. Let $\Under$ be the
\Def{under} operation on $\Reduced \App \SetForests \App \Signature$ such that $\ForestF
\Under \ForestF'$ is obtained by concatenating $\ForestF$ and $\ForestF'$, and then, by
successively grafting the first term of $\ForestF'$ onto the leftmost extremal leaf of the
last term of $\ForestF$, the second term of $\ForestF'$ onto the second leftmost leaf of the
last term $\ForestF$, and so on, until there is no remaining term in $\ForestF'$ or there is
no extremal leaf in the last term of $\ForestF$. Let also $\Over$ be the \Def{over}
operation on $\Reduced \App \SetForests \App \Signature$ defined by $\ForestF \Over
\ForestF' := \ForestF \Conc \ForestF'$.

\begin{Theorem} \label{thm:product_alternative_bases}
    For any positive signature $\Signature$ and any reduced $\Signature$-forests $\ForestF'$
    and $\ForestF''$,
    \begin{equation}
        \BasisE_{\ForestF'} \Product \BasisE_{\ForestF''}
        = \BasisE_{\ForestF' \Over \ForestF''},
        \qquad
        \BasisF_{\ForestF'} \Product \BasisF_{\ForestF''}
        = \sum_{\ForestF \in \Reduced \App \SetForests \App \Signature}
        \Iverson{\ForestF' \Over \ForestF'' \Leq \ForestF \Leq \ForestF' \Under \ForestF''}
        \;
        \BasisF_\ForestF,
        \qquad
        \BasisH_{\ForestF'} \Product \BasisH_{\ForestF''}
        = \BasisH_{\ForestF' \Under \ForestF''}.
    \end{equation}
\end{Theorem}

The product of $\NaturalHopfAlgebra \App \SetTerms \App \Signature$ on the $\BasisF$-basis
is akin to a shuffle of reduced $\Signature$-forests. For instance, in $\NaturalHopfAlgebra
\App \SetTerms \App \SignatureExample$,
\begin{equation}
    \BasisF_{
        \scalebox{.75}{
        \begin{tikzpicture}[Centering,xscale=0.24,yscale=0.22]
            \node[Leaf](0)at(0.00,-2.00){};
            \node[Leaf](2)at(1.00,-2.00){};
            \node[Leaf](3)at(2.00,-2.00){};
            \node[Node](1)at(1.00,0.00){$\GenC$};
            \draw[Edge](0)--(1);
            \draw[Edge](2)--(1);
            \draw[Edge](3)--(1);
            \node[Leaf](r)at(1.00,2){};
            \draw[Edge](r)--(1);
        \end{tikzpicture}}
        \scalebox{.75}{
        \begin{tikzpicture}[Centering,xscale=0.17,yscale=0.26]
            \node[Leaf](0)at(0.00,-1.67){};
            \node[Leaf](2)at(2.00,-3.33){};
            \node[Leaf](4)at(4.00,-3.33){};
            \node[Node](1)at(1.00,0.00){$\GenB$};
            \node[Node](3)at(3.00,-1.67){$\GenB$};
            \draw[Edge](0)--(1);
            \draw[Edge](2)--(3);
            \draw[Edge](3)--(1);
            \draw[Edge](4)--(3);
            \node[Leaf](r)at(1.00,1.67){};
            \draw[Edge](r)--(1);
        \end{tikzpicture}}}
    \Product
    \BasisF_{
        \scalebox{.75}{
        \begin{tikzpicture}[Centering,xscale=0.16,yscale=0.22]
            \node[Leaf](0)at(0.00,-2.00){};
            \node[Leaf](2)at(1.00,-4.00){};
            \node[Leaf](4)at(3.00,-4.00){};
            \node[Leaf](5)at(4.00,-2.00){};
            \node[Node,draw=ColB,fill=ColB!30](1)at(2.00,0.00){$\GenC$};
            \node[Node,draw=ColB,fill=ColB!30](3)at(2.00,-2.00){$\GenA$};
            \draw[Edge](0)--(1);
            \draw[Edge](2)--(3);
            \draw[Edge](3)--(1);
            \draw[Edge](4)--(3);
            \draw[Edge](5)--(1);
            \node[Leaf](r)at(2.00,2){};
            \draw[Edge](r)--(1);
        \end{tikzpicture}}
        \scalebox{.75}{
        \begin{tikzpicture}[Centering,xscale=0.17,yscale=0.4]
            \node[Leaf](1)at(0.00,-1.00){};
            \node[Node,draw=ColB,fill=ColB!30](0)at(0.00,0.00){$\GenA$};
            \draw[Edge](1)--(0);
            \node[Leaf](r)at(0.00,1){};
            \draw[Edge](r)--(0);
        \end{tikzpicture}}
        \scalebox{.75}{
        \begin{tikzpicture}[Centering,xscale=0.17,yscale=0.28]
            \node[Leaf](0)at(0.00,-1.50){};
            \node[Leaf](2)at(2.00,-1.50){};
            \node[Node,draw=ColB,fill=ColB!30](1)at(1.00,0.00){$\GenB$};
            \draw[Edge](0)--(1);
            \draw[Edge](2)--(1);
            \node[Leaf](r)at(1.00,1.5){};
            \draw[Edge](r)--(1);
        \end{tikzpicture}}}
    =
    \BasisF_{
        \scalebox{.75}{
        \begin{tikzpicture}[Centering,xscale=0.24,yscale=0.22]
            \node[Leaf](0)at(0.00,-2.00){};
            \node[Leaf](2)at(1.00,-2.00){};
            \node[Leaf](3)at(2.00,-2.00){};
            \node[Node](1)at(1.00,0.00){$\GenC$};
            \draw[Edge](0)--(1);
            \draw[Edge](2)--(1);
            \draw[Edge](3)--(1);
            \node[Leaf](r)at(1.00,2){};
            \draw[Edge](r)--(1);
        \end{tikzpicture}}
        \scalebox{.75}{
        \begin{tikzpicture}[Centering,xscale=0.17,yscale=0.26]
            \node[Leaf](0)at(0.00,-1.67){};
            \node[Leaf](2)at(2.00,-3.33){};
            \node[Leaf](4)at(4.00,-3.33){};
            \node[Node](1)at(1.00,0.00){$\GenB$};
            \node[Node](3)at(3.00,-1.67){$\GenB$};
            \draw[Edge](0)--(1);
            \draw[Edge](2)--(3);
            \draw[Edge](3)--(1);
            \draw[Edge](4)--(3);
            \node[Leaf](r)at(1.00,1.67){};
            \draw[Edge](r)--(1);
        \end{tikzpicture}}
        \scalebox{.75}{
        \begin{tikzpicture}[Centering,xscale=0.16,yscale=0.22]
            \node[Leaf](0)at(0.00,-2.00){};
            \node[Leaf](2)at(1.00,-4.00){};
            \node[Leaf](4)at(3.00,-4.00){};
            \node[Leaf](5)at(4.00,-2.00){};
            \node[Node,draw=ColB,fill=ColB!30](1)at(2.00,0.00){$\GenC$};
            \node[Node,draw=ColB,fill=ColB!30](3)at(2.00,-2.00){$\GenA$};
            \draw[Edge](0)--(1);
            \draw[Edge](2)--(3);
            \draw[Edge](3)--(1);
            \draw[Edge](4)--(3);
            \draw[Edge](5)--(1);
            \node[Leaf](r)at(2.00,2){};
            \draw[Edge](r)--(1);
        \end{tikzpicture}}
        \scalebox{.75}{
        \begin{tikzpicture}[Centering,xscale=0.17,yscale=0.4]
            \node[Leaf](1)at(0.00,-1.00){};
            \node[Node,draw=ColB,fill=ColB!30](0)at(0.00,0.00){$\GenA$};
            \draw[Edge](1)--(0);
            \node[Leaf](r)at(0.00,1){};
            \draw[Edge](r)--(0);
        \end{tikzpicture}}
        \scalebox{.75}{
        \begin{tikzpicture}[Centering,xscale=0.17,yscale=0.28]
            \node[Leaf](0)at(0.00,-1.50){};
            \node[Leaf](2)at(2.00,-1.50){};
            \node[Node,draw=ColB,fill=ColB!30](1)at(1.00,0.00){$\GenB$};
            \draw[Edge](0)--(1);
            \draw[Edge](2)--(1);
            \node[Leaf](r)at(1.00,1.5){};
            \draw[Edge](r)--(1);
        \end{tikzpicture}}}
    +
    \BasisF_{
        \scalebox{.75}{
        \begin{tikzpicture}[Centering,xscale=0.24,yscale=0.22]
            \node[Leaf](0)at(0.00,-2.00){};
            \node[Leaf](2)at(1.00,-2.00){};
            \node[Leaf](3)at(2.00,-2.00){};
            \node[Node](1)at(1.00,0.00){$\GenC$};
            \draw[Edge](0)--(1);
            \draw[Edge](2)--(1);
            \draw[Edge](3)--(1);
            \node[Leaf](r)at(1.00,2){};
            \draw[Edge](r)--(1);
        \end{tikzpicture}}
        \scalebox{.75}{
        \begin{tikzpicture}[Centering,xscale=0.20,yscale=0.24]
            \node[Leaf](0)at(0.00,-1.80){};
            \node[Leaf](2)at(2.00,-3.60){};
            \node[Leaf](4)at(3.50,-5.40){};
            \node[Leaf](7)at(5.00,-7.20){};
            \node[Leaf](8)at(6.50,-5.40){};
            \node[Node](1)at(1.00,0.00){$\GenB$};
            \node[Node](3)at(3.00,-1.80){$\GenB$};
            \node[Node,draw=ColB,fill=ColB!30](5)at(5.00,-3.60){$\GenC$};
            \node[Node,draw=ColB,fill=ColB!30](6)at(5.00,-5.40){$\GenA$};
            \draw[Edge](0)--(1);
            \draw[Edge](2)--(3);
            \draw[Edge](3)--(1);
            \draw[Edge](4)--(5);
            \draw[Edge](5)--(3);
            \draw[Edge](6)--(5);
            \draw[Edge](7)--(6);
            \draw[Edge](8)--(5);
            \node[Leaf](r)at(1.00,1.8){};
            \draw[Edge](r)--(1);
        \end{tikzpicture}}
        \scalebox{.75}{
        \begin{tikzpicture}[Centering,xscale=0.17,yscale=0.4]
            \node[Leaf](1)at(0.00,-1.00){};
            \node[Node,draw=ColB,fill=ColB!30](0)at(0.00,0.00){$\GenA$};
            \draw[Edge](1)--(0);
            \node[Leaf](r)at(0.00,1){};
            \draw[Edge](r)--(0);
        \end{tikzpicture}}
        \scalebox{.75}{
        \begin{tikzpicture}[Centering,xscale=0.17,yscale=0.28]
            \node[Leaf](0)at(0.00,-1.50){};
            \node[Leaf](2)at(2.00,-1.50){};
            \node[Node,draw=ColB,fill=ColB!30](1)at(1.00,0.00){$\GenB$};
            \draw[Edge](0)--(1);
            \draw[Edge](2)--(1);
            \node[Leaf](r)at(1.00,1.5){};
            \draw[Edge](r)--(1);
        \end{tikzpicture}}}
    +
    \BasisF_{
        \scalebox{.75}{
        \begin{tikzpicture}[Centering,xscale=0.24,yscale=0.22]
            \node[Leaf](0)at(0.00,-2.00){};
            \node[Leaf](2)at(1.00,-2.00){};
            \node[Leaf](3)at(2.00,-2.00){};
            \node[Node](1)at(1.00,0.00){$\GenC$};
            \draw[Edge](0)--(1);
            \draw[Edge](2)--(1);
            \draw[Edge](3)--(1);
            \node[Leaf](r)at(1.00,2){};
            \draw[Edge](r)--(1);
        \end{tikzpicture}}
        \scalebox{.75}{
        \begin{tikzpicture}[Centering,xscale=0.22,yscale=0.24]
            \node[Leaf](0)at(1.00,-1.80){};
            \node[Leaf](2)at(1.50,-5.40){};
            \node[Leaf](5)at(3.00,-7.20){};
            \node[Leaf](6)at(4.50,-5.40){};
            \node[Leaf](8)at(6.00,-3.60){};
            \node[Node](1)at(2.00,0.00){$\GenB$};
            \node[Node,draw=ColB,fill=ColB!30](3)at(3.00,-3.60){$\GenC$};
            \node[Node,draw=ColB,fill=ColB!30](4)at(3.00,-5.40){$\GenA$};
            \node[Node](7)at(5.00,-1.80){$\GenB$};
            \draw[Edge](0)--(1);
            \draw[Edge](2)--(3);
            \draw[Edge](3)--(7);
            \draw[Edge](4)--(3);
            \draw[Edge](5)--(4);
            \draw[Edge](6)--(3);
            \draw[Edge](7)--(1);
            \draw[Edge](8)--(7);
            \node[Leaf](r)at(2.00,1.8){};
            \draw[Edge](r)--(1);
        \end{tikzpicture}}
        \scalebox{.75}{
        \begin{tikzpicture}[Centering,xscale=0.17,yscale=0.4]
            \node[Leaf](1)at(0.00,-1.00){};
            \node[Node,draw=ColB,fill=ColB!30](0)at(0.00,0.00){$\GenA$};
            \draw[Edge](1)--(0);
            \node[Leaf](r)at(0.00,1){};
            \draw[Edge](r)--(0);
        \end{tikzpicture}}
        \scalebox{.75}{
        \begin{tikzpicture}[Centering,xscale=0.17,yscale=0.28]
            \node[Leaf](0)at(0.00,-1.50){};
            \node[Leaf](2)at(2.00,-1.50){};
            \node[Node,draw=ColB,fill=ColB!30](1)at(1.00,0.00){$\GenB$};
            \draw[Edge](0)--(1);
            \draw[Edge](2)--(1);
            \node[Leaf](r)at(1.00,1.5){};
            \draw[Edge](r)--(1);
        \end{tikzpicture}}}
    +
    \BasisF_{
        \scalebox{.75}{
        \begin{tikzpicture}[Centering,xscale=0.24,yscale=0.22]
            \node[Leaf](0)at(0.00,-2.00){};
            \node[Leaf](2)at(1.00,-2.00){};
            \node[Leaf](3)at(2.00,-2.00){};
            \node[Node](1)at(1.00,0.00){$\GenC$};
            \draw[Edge](0)--(1);
            \draw[Edge](2)--(1);
            \draw[Edge](3)--(1);
            \node[Leaf](r)at(1.00,2){};
            \draw[Edge](r)--(1);
        \end{tikzpicture}}
        \scalebox{.75}{
        \begin{tikzpicture}[Centering,xscale=0.22,yscale=0.23]
            \node[Leaf](0)at(1.00,-2.00){};
            \node[Leaf](2)at(1.50,-6.00){};
            \node[Leaf](5)at(3.00,-8.00){};
            \node[Leaf](6)at(4.50,-6.00){};
            \node[Leaf](9)at(6.00,-6.00){};
            \node[Node](1)at(2.00,0.00){$\GenB$};
            \node[Node,draw=ColB,fill=ColB!30](3)at(3.00,-4.00){$\GenC$};
            \node[Node,draw=ColB,fill=ColB!30](4)at(3.00,-6.00){$\GenA$};
            \node[Node](7)at(5.00,-2.00){$\GenB$};
            \node[Node,draw=ColB,fill=ColB!30](8)at(6.00,-4.00){$\GenA$};
            \draw[Edge](0)--(1);
            \draw[Edge](2)--(3);
            \draw[Edge](3)--(7);
            \draw[Edge](4)--(3);
            \draw[Edge](5)--(4);
            \draw[Edge](6)--(3);
            \draw[Edge](7)--(1);
            \draw[Edge](8)--(7);
            \draw[Edge](9)--(8);
            \node[Leaf](r)at(2.00,2){};
            \draw[Edge](r)--(1);
        \end{tikzpicture}}
        \scalebox{.75}{
        \begin{tikzpicture}[Centering,xscale=0.17,yscale=0.28]
            \node[Leaf](0)at(0.00,-1.50){};
            \node[Leaf](2)at(2.00,-1.50){};
            \node[Node,draw=ColB,fill=ColB!30](1)at(1.00,0.00){$\GenB$};
            \draw[Edge](0)--(1);
            \draw[Edge](2)--(1);
            \node[Leaf](r)at(1.00,1.5){};
            \draw[Edge](r)--(1);
        \end{tikzpicture}}}.
\end{equation}

These three bases, together with the over and under operations, and the partial order on
reduced $\Signature$-forests mimic classical basis constructions of Hopf
algebras~\cite{DHT02,LR02,HNT05}.

\section{Quotient operads and Hopf subalgebras} \label{sec:hopf_subalgebras}
Let $\Signature$ be a signature and $\Equiv$ be an operad congruence of the free operad
$\SetTerms \App \Signature$. The equivalence relation $\Equiv$ is extended on $\Reduced \App
\SetForests \App \Signature$ by setting $\ForestF \Equiv \ForestF'$ if $\Length \App
\ForestF = \Length \App \ForestF'$ and, for any $i \in [\Length \App \ForestF]$, $\ForestF
\App i \Equiv \ForestF' \App i$. The operad congruence $\Equiv$ is \Def{compatible with the
degree} if $\TermT \Equiv \TermT'$ implies $\Degree \App \TermT = \Degree \App \TermT'$ for
any $\TermT, \TermT' \in \SetTerms \App \Signature$. The operad congruence $\Equiv$ is of
\Def{finite type} if for any $\TermT \in \SetTerms \App \Signature$, the
$\Equiv$-equivalence class $[\TermT]_{\Equiv}$ of $\TermT$ is finite. In this case, for any
$\ForestF \in \Reduced \App \SetForests \App \Signature$,
\begin{equation}
    \BasisE_{[\ForestF]_{\Equiv}} :=
    \sum_{\ForestF' \in [\ForestF]_{\Equiv}} \BasisE_{\ForestF'}
\end{equation}
is a well-defined element of $\NaturalHopfAlgebra \App \SetTerms \App \Signature$.

\begin{Theorem} \label{thm:quotient_operad_hopf_subalgebra}
    Let $\Signature$ be a signature and $\Equiv$ be an operad congruence of $\SetTerms \App
    \Signature$ compatible with the degree and of finite type. By denoting by $\Space$ the
    linear span of the $\BasisE_{[\ForestF]_{\Equiv}}$, $\ForestF \in \Reduced \App
    \SetForests \App \Signature$, the following assertions hold:
    \begin{enumerate}[label=({\sf \roman*})]
        \item \label{item:quotient_operad_hopf_subalgebra_1}
        The space $\Space$ is a Hopf sub-algebra of $\NaturalHopfAlgebra \App \SetTerms \App
        \Signature$.
        \item \label{item:quotient_operad_hopf_subalgebra_2}
        The Hopf algebra $\Space$ is isomorphic to $\NaturalHopfAlgebra \App \Par{\SetTerms
        \App \Signature} /_{\Equiv}$.
    \end{enumerate}
\end{Theorem}

Point~\ref{item:quotient_operad_hopf_subalgebra_1} of
Theorem~\ref{thm:quotient_operad_hopf_subalgebra} is a particular case of~\cite[Theorem
3.13]{BG16}. Point~\ref{item:quotient_operad_hopf_subalgebra_2} provides an interpretation
of $\NaturalHopfAlgebra \App \Par{\SetTerms \App \Signature} /_{\Equiv}$ as a Hopf
subalgebra of the natural Hopf algebra of a free operad. The condition of compatibility with
the degree for $\Equiv$ is important for the fact that $\NaturalHopfAlgebra \App
\Par{\SetTerms \App \Signature} /_{\Equiv}$ inherits the graduation of $\NaturalHopfAlgebra
\App \SetTerms \App \Signature$.

\section{Polynomial realization} \label{sec:polynomial_realization}
Given a signature $\Signature$, the class of \Def{$\Signature$-forest-like alphabets} is the
class of related alphabets $A$ endowed with relations $\RootRelation$,
$\DecorationRelation_\GenG$, and $\EdgeRelation{j}$ such that
\begin{enumerate}[label=({\sf \roman*})]
    \item $\RootRelation$ is a unary relation called \Def{root relation};
    \item for any $\GenG \in \Signature$, $\DecorationRelation_\GenG$ is a unary relation
    called \Def{$\GenG$-decoration relation};
    \item for any $j \geq 1$, $\EdgeRelation{j}$ is a binary relation called \Def{$j$-edge
    relation}.
\end{enumerate}

Let $A'$ and $A''$ be two $\Signature$-forest-like alphabets, each endowed respectively with
the relations $\RootRelation'$, $\DecorationRelation'_\GenG$, and $\EdgeRelation{j}'$, and
$\RootRelation''$, $\DecorationRelation''_\GenG$, and $\EdgeRelation{j}''$. The
\Def{disjoint sum} of $A'$ and $A''$ is the $\Signature$-forest-like alphabet $A := A'
\AlphabetSum A''$ endowed with the relations $\RootRelation$, $\DecorationRelation_\GenG$,
and $\EdgeRelation{j}$ such that $A := A' \sqcup A''$, $\RootRelation := \RootRelation'
\sqcup \RootRelation''$, $\DecorationRelation_\GenG := \DecorationRelation_\GenG' \sqcup
\DecorationRelation_\GenG''$, and $a_1 \EdgeRelation{j} a_2$ holds if $a_1, a_2 \in A'$ and
$a_1 \EdgeRelation{j}' a_2$, or $a_1, a_2 \in A''$ and $a_1 \EdgeRelation{j}'' a_2$, or $a_1
\in A'$, $a_2 \in A''$, and $a_2 \in \RootRelation''$.

Let $A$ be an $\Signature$-forestlike alphabet. Given $\ForestF \in \Reduced \App
\SetForests \App \Signature$, $w \in A^*$ is \Def{$\ForestF$-compatible} if the following
four assertions are satisfied:
\begin{enumerate*}[label=({\sf \roman*})]
    \item $\Length \App w = \Degree \App \ForestF$;
    \item for any $i \in [\Degree \App \ForestF]$, $i$ is a root of $\ForestF$ implies $w
    \App i \in \RootRelation$;
    \item for any $i \in [\Degree \App \ForestF]$ the internal node $i$ of $\ForestF$ is
    decorated by $\GenG \in \Signature$ implies $w \App i \in \DecorationRelation_\GenG$;
    \item for any $i, i' \in [\Degree \App \ForestF]$, $i \Edge{\ForestF}{j} i'$ implies $w
    \App i \EdgeRelation{j} w \App i'$.
\end{enumerate*}
This property is denoted by $w \CompatibleWord^A \ForestF$. For instance, by considering the
reduced $\SignatureExample$-forest $\ForestF$ of Figure~\ref{fig:example_forest}, any
$\ForestF$-compatible word $w \in A^*$ satisfies $\Length \App w = 7$, $w \App 1, w \App 5
\in \RootRelation$, $w \App 2, w \App 4, w \App 7 \in \DecorationRelation_\GenA$, $w \App 3,
w \App 5, w \App 6 \in \DecorationRelation_\GenB$, $w \App 1 \in \DecorationRelation_\GenC$,
$w \App 1 \EdgeRelation{1} w \App 2$, $w \App 1 \EdgeRelation{3} w \App 3$, $w \App 3
\EdgeRelation{1} w \App 4$, $w \App 5 \EdgeRelation{2} w \App 6$, and $w \App 6
\EdgeRelation{1} w \App 7$.

Let $\Realization_A : \NaturalHopfAlgebra \App \SetTerms \App \Signature \to \K \Angle{A}$
be the linear map defined for any $\ForestF \in \Reduced \App \SetForests \App \Signature$
by
\begin{equation} \label{equ:realization_forest}
    \Realization_A \App \BasisE_{\ForestF}
    := \sum_{w \in A^*} \Iverson{w \CompatibleWord^A \ForestF} \; w.
\end{equation}
The $A$-polynomial $\Realization_A \App \BasisE_\ForestF$ is the \Def{$A$-realization} of
$\ForestF$.

We now construct a specific $\Signature$-forest-like alphabet $\Alphabet_\Signature$ in
order to obtain an injective map $\Realization_{\Alphabet_\Signature}$. For any signature
$\Signature$, let the $\Signature$-forest-like alphabet
\begin{math}
    \Alphabet_\Signature
\end{math}
such that
\begin{equation}
    \Alphabet_\Signature
    := \Bra{\Letter_{\GenG, u} : \GenG \in \Signature \mbox{ and } u \in \N^*},
\end{equation}
the root relation is defined by
\begin{math}
    \RootRelation := \Bra{\Letter_{\GenG, u} \in \Alphabet_\Signature :
    u = 0^\ell, \ell \geq 0},
\end{math}
the $\GenG$-decoration relation $\DecorationRelation_\GenG$ is defined by
\begin{math}
    \DecorationRelation_\GenG :=
    \Bra{\Letter_{\GenG', u} \in \Alphabet_\Signature : \GenG' = \GenG},
\end{math}
the $j$-edge relation $\EdgeRelation{j}$ is defined by
\begin{math}
    \Letter_{\GenG, u} \EdgeRelation{j} \Letter_{\GenG', u \Conc j \Conc 0^\ell}
\end{math}
for any $\GenG, \GenG' \in \Signature$, $u \in \N^*$, and $\ell \geq 0$. Observe that these
binary relations are antisymmetric but are neither transitive nor reflexive. For instance,
by considering the reduced $\SignatureExample$-forest $\ForestF$ of
Figure~\ref{fig:example_forest},
\begin{equation}
    \Realization_{\Alphabet_\Signature} \App \BasisE_\ForestF
    =
    \sum_{\ell_1, \dots, \ell_7 \in \N}
    \Letter_{\GenC, 0^{\ell_1}}
    \;
    \Letter_{\GenA, 0^{\ell_1} 1 0^{\ell_2}}
    \;
    \Letter_{\GenB, 0^{\ell_1} 3 0^{\ell_3}}
    \;
    \Letter_{\GenA, 0^{\ell_1} 3 0^{\ell_3} 1 0^{\ell_4}}
    \;
    \Letter_{\GenB, 0^{\ell_5}}
    \;
    \Letter_{\GenB, 0^{\ell_5} 2 0^{\ell_6}}
    \;
    \Letter_{\GenB, 0^{\ell_5} 2 0^{\ell_6} 1 0^{\ell_7}}.
\end{equation}

\begin{Theorem} \label{thm:realization}
    For any signature $\Signature$, the class of $\Signature$-forest-like alphabets,
    together with the alphabet disjoint sum operation $\AlphabetSum$, the map
    $\Realization_A$, and the alphabet $\Alphabet_\Signature$, forms a polynomial
    realization of the Hopf algebra $\NaturalHopfAlgebra \App \SetTerms \App \Signature$.
\end{Theorem}

By Theorems~\ref{thm:quotient_operad_hopf_subalgebra} and~\ref{thm:realization}, the map
$\Realization_{\Alphabet_\Signature}$ provides also a polynomial realization of the Hopf
subalgebra $\NaturalHopfAlgebra \App \Par{\SetTerms \App \Signature}/_{\Equiv}$ of
$\NaturalHopfAlgebra \App \SetTerms \App \Signature$ when $\Equiv$ is an operad congruence
of $\SetTerms \App \Signature$ compatible with the degree and of finite type.

\section{Word quasi-symmetric functions Hopf algebra} \label{sec:wqsym}
Let $A$ be a totally ordered alphabet. The \Def{packing} of $w \in A^*$ is the word $\Pack
\App w$ on positive integers obtained by replacing each letter $a$ of $w$ by the number of
distinct letters less than or equal to $a$ in $w$. For instance, for $A = \N$, $\Pack \App
4234473 = 3123342$. When $w$ is a word of positive integers such that $\Pack \App w = w$,
$w$ is \Def{packed}. The set of packed words is denoted by $\SetPackedWords$. The Hopf
algebra of \Def{word quasi-symmetric functions} $\WQSym$~\cite{Hiv99,NT06} is defined on the
linear span of $\SetPackedWords$. The \Def{fundamental basis} of $\WQSym$ is the set
$\Bra{\BasisM_u : u \in \SetPackedWords}$. We need here only to describe the polynomial
realization of $\WQSym$, which satisfies
\begin{equation}
    \Realization_A \App \BasisM_u = \sum_{w \in A^*} \Iverson{\Pack(w) = u} \; w.
\end{equation}
By setting  $\Alphabet$ as the related alphabet $\Bra{\Letter_i : i \in \N}$ totally ordered
by the order relation $\Leq$ satisfying $\Letter_i \Leq \Letter_{i'}$ if $i \leq i'$, we
have for instance
\begin{equation}
    \Realization_\Alphabet \App \BasisM_{2113}
    = \sum_{\ell_1, \ell_2, \ell_3, \ell_4 \in \N}
    \Iverson{\ell_2 = \ell_3 < \ell_1 < \ell_4} \;
    \Letter_{\ell_1} \Letter_{\ell_2} \Letter_{\ell_3} \Letter_{\ell_4}.
\end{equation}

Let $\AlphabetN$ be the $\Signature$-forest-like alphabet such that $\AlphabetN :=
\Alphabet$, $\RootRelation := \AlphabetN$, $\DecorationRelation_\GenG := \AlphabetN$, and
$\Letter_i \EdgeRelation{j} \Letter_{i'}$ if $i < i'$. Given $\ForestF \in \Reduced \App
\SetForests \App \Signature$, $u \in \SetPackedWords$ is \Def{$\ForestF$-compatible} if
$\Letter_{u \App 1} \dots \Letter_{u \App \Length \App u} \CompatibleWord^{\AlphabetN}
\ForestF$. This property is denoted by $u \CompatiblePackedWord \ForestF$. For instance, by
considering the reduced $\SignatureExample$-forest $\ForestF$ defined in
Figure~\ref{fig:example_forest}, $1223123 \CompatiblePackedWord \ForestF$ and $1525346
\CompatiblePackedWord \ForestF$.

\begin{Theorem} \label{thm:map_to_wqsym}
    For any signature $\Signature$, the space
    \begin{math}
        \Realization_{\AlphabetN} \App \NaturalHopfAlgebra \App \SetTerms \App \Signature
    \end{math}
    is a Hopf subalgebra of $\Realization_\Alphabet \App \WQSym$. More precisely, for any
    $\ForestF \in \Reduced \App \SetForests \App \Signature$,
    \begin{equation}
        \Realization_{\AlphabetN} \App \BasisE_\ForestF
        =
        \sum_{u \in \SetPackedWords}
        \Iverson{u \CompatiblePackedWord \ForestF} \;
        \Realization_\Alphabet \App \BasisM_u.
    \end{equation}
\end{Theorem}

By Theorems~\ref{thm:quotient_operad_hopf_subalgebra} and~\ref{thm:map_to_wqsym},
$\Realization_{\Alphabet_\N} \App \NaturalHopfAlgebra \App \Par{\SetTerms \App
\Signature}/_{\Equiv}$ is a Hopf subalgebra of $\Realization_\Alphabet \App \WQSym$ when
$\Equiv$ is an operad congruence compatible with the degree and of finite type. Observe that
the map $\Realization_{\AlphabetN}$ is not injective as, for instance,
\begin{equation}
    \Realization_{\AlphabetN} \App
    \BasisE_{
        \scalebox{.75}{
        \begin{tikzpicture}[Centering,xscale=0.26,yscale=0.21]
            \node[Leaf](1)at(0.00,-4.00){};
            \node[Leaf](3)at(2.00,-4.00){};
            \node[Leaf](5)at(4.00,-4.00){};
            \node[Node](0)at(0.00,-2.00){$\GenA$};
            \node[Node](2)at(1.50,0.00){$\GenB$};
            \node[Node](4)at(3.00,-2.00){$\GenB$};
            \draw[Edge](0)--(2);
            \draw[Edge](1)--(0);
            \draw[Edge](3)--(4);
            \draw[Edge](4)--(2);
            \draw[Edge](5)--(4);
            \node[Leaf](r)at(1.50,2){};
            \draw[Edge](r)--(2);
        \end{tikzpicture}}}
    =
    \BasisM_{122} + \BasisM_{123} + \BasisM_{132}
    =
    \Realization_{\AlphabetN} \App
    \BasisE_{
        \scalebox{.75}{
        \begin{tikzpicture}[Centering,xscale=0.24,yscale=0.17]
            \node[Leaf](0)at(-1.00,-2.67){};
            \node[Leaf](2)at(1.00,-5.33){};
            \node[Leaf](4)at(3.00,-5.33){};
            \node[Leaf](5)at(4.00,-5.33){};
            \node[Leaf](7)at(6.00,-5.33){};
            \node[Node](1)at(2.00,0.00){$\GenC$};
            \node[Node](3)at(2.00,-2.67){$\GenB$};
            \node[Node](6)at(5.00,-2.67){$\GenB$};
            \draw[Edge](0)--(1);
            \draw[Edge](2)--(3);
            \draw[Edge](3)--(1);
            \draw[Edge](4)--(3);
            \draw[Edge](5)--(6);
            \draw[Edge](6)--(1);
            \draw[Edge](7)--(6);
            \node[Leaf](r)at(2.00,2.67){};
            \draw[Edge](r)--(1);
        \end{tikzpicture}}}.
\end{equation}

\section{Connes-Kreimer and decorated forests Hopf algebras} \label{sec:connes_kreimer}
For any set $D$, a \Def{$D$-decorated forest} is a word of nonempty planar rooted trees
where nodes are decorated on $D$. In~\cite{Foi02a}, Foissy introduces a Hopf algebra
structure $\NCK_D$ on the linear span of $D$-decorated forests. This is a generalization of
the noncommutative Connes-Kreimer Hopf algebra as the original case occurs when $D$ is a
singleton.

For any set $D$, let the signature $\Signature_D := \Bra{\alpha_{d, n} : d \in D, n \in \N}$
such that $\Arity \App \alpha_{d, n} = n$. Let also $\Alphabet_{D, \N}$ be the
$\Signature_D$-forest-like alphabet such that $\Alphabet_{D, \N} := \Bra{\Letter_{d, i} : d
\in D, i \in \N}$, $\RootRelation := \Alphabet_{D, \N}$, $\DecorationRelation_{\alpha_{d,
n}} := \Bra{\Letter_{d, i} \in \Alphabet_{D, \N} : i \in \N}$, and $\Letter_{d, i}
\EdgeRelation{j} \Letter_{d', i'}$ if $i < i'$.

\begin{Theorem} \label{thm:map_to_nck}
    For any set $D$, the Hopf algebra $\Realization_{\Alphabet_{D, \N}} \App
    \NaturalHopfAlgebra \App \SetTerms \App \Signature_D$ is isomorphic to the Hopf algebra
    $\NCK_D$ of $D$-decorated forests.
\end{Theorem}

Theorem~\ref{thm:map_to_nck} provides a polynomial realization of $\NCK_D$. Indeed, given a
$D$-decorated forest $F$, the $\Alphabet_{D, \N}$-polynomial associated with the basis
element $\BasisE_F$ of $\NCK_D$, is $\Realization_{\Alphabet_{D, \N}} \App \BasisE_\ForestF$
where $\ForestF$ is the $\Signature_D$-forest obtained by replacing each internal node
decorated by $d \in D$ and having $n$ children by $\alpha_{d, n}$. For instance, for $D :=
\Bra{1, 2, 3}$,
\begin{align}
    \Realization_{\Alphabet_{D, \N}}
    & \App
    \BasisE_{
        \scalebox{.75}{
        \begin{tikzpicture}[Centering,xscale=0.32,yscale=0.28]
            \node[Node,draw=ColA!50!ColB,fill=ColA!25!ColB!10](0)at(0.00,-2.00){$2$};
            \node[Node,draw=ColA!50!ColB,fill=ColA!25!ColB!10](2)at(1.50,-4.00){$2$};
            \node[Node,draw=ColA!50!ColB,fill=ColA!25!ColB!10](4)at(3.00,-4.00){$3$};
            \node[Node,draw=ColA!50!ColB,fill=ColA!25!ColB!10](5)at(4.50,-4.00){$1$};
            \node[Node,draw=ColA!50!ColB,fill=ColA!25!ColB!10](1)at(1.50,0.00){$3$};
            \node[Node,draw=ColA!50!ColB,fill=ColA!25!ColB!10](3)at(3.00,-2.00){$3$};
            \draw[Edge](0)--(1);
            \draw[Edge](2)--(3);
            \draw[Edge](3)--(1);
            \draw[Edge](4)--(3);
            \draw[Edge](5)--(3);
        \end{tikzpicture}}
        \;
        \scalebox{.75}{
        \begin{tikzpicture}[Centering,xscale=0.17,yscale=0.5]
            \node[Node,draw=ColA!50!ColB,fill=ColA!25!ColB!10](1)at(0.00,-1.00){$2$};
            \node[Node,draw=ColA!50!ColB,fill=ColA!25!ColB!10](0)at(0.00,0.00){$1$};
            \draw[Edge](1)--(0);
        \end{tikzpicture}}}
    :=
    \Realization_{\Alphabet_{D, \N}}
    \App
    \BasisE_{
        \scalebox{.75}{
        \begin{tikzpicture}[Centering,xscale=0.5,yscale=0.34]
            \node[Node](0)at(0.00,-2.00){$\alpha_{2, 0}$};
            \node[Node](2)at(1.25,-4.00){$\alpha_{2, 0}$};
            \node[Node](4)at(3.00,-4.00){$\alpha_{3, 0}$};
            \node[Node](5)at(4.75,-4.00){$\alpha_{1, 0}$};
            \node[Node](1)at(1.50,0.00){$\alpha_{3, 2}$};
            \node[Node](3)at(3.00,-2.00){$\alpha_{3, 2}$};
            \draw[Edge](0)--(1);
            \draw[Edge](2)--(3);
            \draw[Edge](3)--(1);
            \draw[Edge](4)--(3);
            \draw[Edge](5)--(3);
            \node[Leaf](r)at(1.50,2){};
            \draw[Edge](r)--(1);
        \end{tikzpicture}}
        \;
        \scalebox{.75}{
        \begin{tikzpicture}[Centering,xscale=0.17,yscale=0.8]
            \node[Node](1)at(0.00,-1.00){$\alpha_{2, 0}$};
            \node[Node](0)at(0.00,0.00){$\alpha_{1, 1}$};
            \draw[Edge](1)--(0);
            \node[Leaf](r)at(0.00,0.75){};
            \draw[Edge](r)--(0);
        \end{tikzpicture}}}
    \\
    & =
    \sum_{\ell_1, \dots, \ell_8 \in \N}
    \Iverson{\ell_1 < \ell_2, \ell_3}
    \Iverson{\ell_3 < \ell_4, \ell_5, \ell_6}
    \Iverson{\ell_7 < \ell_8} \;
    \Letter_{3, \ell_1}
    \Letter_{2, \ell_2}
    \Letter_{3, \ell_3}
    \Letter_{2, \ell_4}
    \Letter_{3, \ell_5}
    \Letter_{1, \ell_6}
    \Letter_{1, \ell_7}
    \Letter_{2, \ell_8}.
    \notag
\end{align}

\section{Faà di Bruno Hopf algebras} \label{sec:faa_di_bruno}
Let $\MAs_\Signature$ be the quotient of the free operad $\SetTerms \App \Signature$ by the
congruence $\Equiv$ satisfying, for any $\GenG, \GenG' \in \Signature$, $i \in [\Arity \App
\GenG]$, and $i' \in \Han{\Arity \App \GenG'}$, $\GenG \circ_i \GenG' \Equiv \GenG'
\circ_{i'} \GenG$. These operads $\MAs_\Signature$ are generalizations of multiassociative
operads, introduced in~\cite{Gir16}, which contain originally only binary generators. For
this reason, we call $\MAs_\Signature$ the \Def{$\Signature$-multi-multiassociative operad}.

\begin{Proposition} \label{prop:realization_operad_mas}
    For any signature $\Signature$, the operad $\MAs_\Signature$ admits the following
    realization. For any $n \geq 0$, $\MAs_\Signature \App n$ is the set of multisets
    $\Bag{\GenG_1, \dots, \GenG_\ell}$ on $\Signature$ such that $\Arity \App \GenG_1 +
    \dots + \Arity \App \GenG_\ell - \ell + 1 = n$. Moreover, for any $\MultisetM,
    \MultisetM' \in \MAs_\Signature$ and $i \in [\Arity \App \MultisetM]$, $\MultisetM
    \circ_i \MultisetM'$ is the union of $\MultisetM$ and $\MultisetM'$.
\end{Proposition}

Recall that $\NCFdB$ is the noncommutative Faà di Bruno Hopf algebra whose construction is
detailed in Section~\ref{subsec:natural_hopf_algebras}. In~\cite{Foi08}, a deformation
$\NCFdB_r$ of $\NCFdB$ is introduced. The special case $\NCFdB_0$ is the Hopf algebra of
noncommutative symmetric functions $\NCSF$, and for any $r \geq 1$, each $\NCFdB_r$ is
isomorphic (in a nontrivial way) to $\NCFdB$.

\begin{Proposition} \label{prop:natural_hopf_algebra_mas_ncfdb}
    Let $\Signature$ be a signature of profile $0^{r + 1} 1 0 \dots$ where $r \in \N$. The
    Hopf algebra $\NaturalHopfAlgebra \App \MAs_{\Signature}$ is isomorphic to the
    noncommutative Faà di Bruno Hopf algebra $\NCFdB_r$.
\end{Proposition}

This result is a consequence of an alternative construction of $\NCFdB_r$ via pros, whose
details can be found in~\cite{BG16}.

\section{Multi-symmetric functions Hopf algebras} \label{sec:multi_symmetric_functions}
Let us now consider the Hopf algebra $\NCSF^{(s)}$ of noncommutative multi-symmetric
functions of level $s \in \N$. This Hopf algebra is a noncommutative analog of
multi-symmetric functions, introduced in~\cite{NT10}.

\begin{Proposition} \label{prop:natural_hopf_algebra_mas_multi_sym}
    Let $\Signature$ be a signature of profile $0 s 0 \dots$ where $s \in \N$. The Hopf
    algebra $\NaturalHopfAlgebra\App \MAs_{\Signature}$ is isomorphic to the noncommutative
    multi-symmetric functions Hopf algebra $\NCSF^{(s)}$.
\end{Proposition}

This framework offered by the multi-multiassociative operads leads naturally to Hopf
algebras which are to $\NCSF^{(s)}$ what $\NCFdB_r$ is to $\NCFdB_0$. These Hopf algebras
depend on the two parameters $r, s \in \N$, are denoted by $\NCFdB_{r, s}$, and are defined
as $\NaturalHopfAlgebra \App \MAs_{\Signature}$ where the $\Signature$ are signatures of
profiles $0^{r + 1} s 0 \dots$. By Theorems~\ref{thm:quotient_operad_hopf_subalgebra}
and~\ref{thm:realization}, we have, for free, polynomial realizations of $\NCFdB_{r, s}$ on
the alphabet $\Alphabet_\Signature$.

By the following result, these polynomial realizations can expressed on a simpler alphabet
$\Alphabet_{\Signature, \N}$ defined as follows. Let $\Alphabet_{\Signature, \N}$ be the
$\Signature$-forest-like alphabet such that $\Alphabet_{\Signature, \N} :=
\Bra{\Letter_{\GenG, i} : \GenG \in \Signature, i \in \N}$, $\RootRelation :=
\Alphabet_{\Signature, \N}$, $\DecorationRelation_\GenG := \Bra{\Letter_{\GenG, i} \in
\Alphabet_{\Signature, \N}: i \in \N}$, and $\Letter_{\GenG, i} \EdgeRelation{j}
\Letter_{\GenG', i'}$ if $i < i'$.

\begin{Theorem} \label{thm:injective_length_morphism_fdb}
    Let $\Signature$ be a signature of profile $0^{r + 1} s 0 \dots$ where $r, s \in \N$.
    The Hopf algebra $\Realization_{\Alphabet_{\Signature, \N}} \App \NaturalHopfAlgebra
    \App \MAs_\Signature$ is isomorphic to $\NCFdB_{r, s}$.
\end{Theorem}

For instance, let $\Signature := \Bra{\GenA, \GenB}$ be a signature of profile $0^{r + 1} s
0 \dots$ where $r = 1$ and $s = 2$. In $\NCFdB_{r, s}$, we have
\begin{align}
    \Realization_{\Alphabet_{\Signature, \N}} \App
    \BasisE_{\Bag{\GenA \GenA \GenB} \Bag{\GenA}}
    & =
    4 \; \Realization_{\Alphabet_{\Signature, \N}} \App
    \BasisE_{
        \scalebox{.75}{
        \begin{tikzpicture}[Centering,xscale=0.17,yscale=0.24]
            \node[Leaf](0)at(0.00,-5.25){};
            \node[Leaf](2)at(2.00,-5.25){};
            \node[Leaf](4)at(4.00,-3.50){};
            \node[Leaf](6)at(6.00,-1.75){};
            \node[Node](1)at(1.00,-3.50){$\GenB$};
            \node[Node](3)at(3.00,-1.75){$\GenA$};
            \node[Node](5)at(5.00,0.00){$\GenA$};
            \draw[Edge](0)--(1);
            \draw[Edge](1)--(3);
            \draw[Edge](2)--(1);
            \draw[Edge](3)--(5);
            \draw[Edge](4)--(3);
            \draw[Edge](6)--(5);
            \node[Leaf](r)at(5.00,1.75){};
            \draw[Edge](r)--(5);
        \end{tikzpicture}}
        \scalebox{.75}{
        \begin{tikzpicture}[Centering,xscale=0.17,yscale=0.24]
            \node[Leaf](0)at(0.00,-1.50){};
            \node[Leaf](2)at(2.00,-1.50){};
            \node[Node](1)at(1.00,0.00){$\GenA$};
            \draw[Edge](0)--(1);
            \draw[Edge](2)--(1);
            \node[Leaf](r)at(1.00,1.5){};
            \draw[Edge](r)--(1);
        \end{tikzpicture}}}
    +
    4 \; \Realization_{\Alphabet_{\Signature, \N}} \App
    \BasisE_{
        \scalebox{.75}{
        \begin{tikzpicture}[Centering,xscale=0.17,yscale=0.24]
            \node[Leaf](0)at(0.00,-5.25){};
            \node[Leaf](2)at(2.00,-5.25){};
            \node[Leaf](4)at(4.00,-3.50){};
            \node[Leaf](6)at(6.00,-1.75){};
            \node[Node](1)at(1.00,-3.50){$\GenA$};
            \node[Node](3)at(3.00,-1.75){$\GenB$};
            \node[Node](5)at(5.00,0.00){$\GenA$};
            \draw[Edge](0)--(1);
            \draw[Edge](1)--(3);
            \draw[Edge](2)--(1);
            \draw[Edge](3)--(5);
            \draw[Edge](4)--(3);
            \draw[Edge](6)--(5);
            \node[Leaf](r)at(5.00,1.75){};
            \draw[Edge](r)--(5);
        \end{tikzpicture}}
        \scalebox{.75}{
        \begin{tikzpicture}[Centering,xscale=0.17,yscale=0.24]
            \node[Leaf](0)at(0.00,-1.50){};
            \node[Leaf](2)at(2.00,-1.50){};
            \node[Node](1)at(1.00,0.00){$\GenA$};
            \draw[Edge](0)--(1);
            \draw[Edge](2)--(1);
            \node[Leaf](r)at(1.00,1.5){};
            \draw[Edge](r)--(1);
        \end{tikzpicture}}}
    +
    4 \; \Realization_{\Alphabet_{\Signature, \N}} \App
    \BasisE_{
        \scalebox{.75}{
        \begin{tikzpicture}[Centering,xscale=0.17,yscale=0.24]
            \node[Leaf](0)at(0.00,-5.25){};
            \node[Leaf](2)at(2.00,-5.25){};
            \node[Leaf](4)at(4.00,-3.50){};
            \node[Leaf](6)at(6.00,-1.75){};
            \node[Node](1)at(1.00,-3.50){$\GenA$};
            \node[Node](3)at(3.00,-1.75){$\GenA$};
            \node[Node](5)at(5.00,0.00){$\GenB$};
            \draw[Edge](0)--(1);
            \draw[Edge](1)--(3);
            \draw[Edge](2)--(1);
            \draw[Edge](3)--(5);
            \draw[Edge](4)--(3);
            \draw[Edge](6)--(5);
            \node[Leaf](r)at(5.00,1.75){};
            \draw[Edge](r)--(5);
        \end{tikzpicture}}
        \scalebox{.75}{
        \begin{tikzpicture}[Centering,xscale=0.17,yscale=0.24]
            \node[Leaf](0)at(0.00,-1.50){};
            \node[Leaf](2)at(2.00,-1.50){};
            \node[Node](1)at(1.00,0.00){$\GenA$};
            \draw[Edge](0)--(1);
            \draw[Edge](2)--(1);
            \node[Leaf](r)at(1.00,1.5){};
            \draw[Edge](r)--(1);
        \end{tikzpicture}}}
    \\
    & \quad +
    \Realization_{\Alphabet_{\Signature, \N}} \App
    \BasisE_{
        \scalebox{.75}{
        \begin{tikzpicture}[Centering,xscale=0.17,yscale=0.18]
            \node[Leaf](0)at(0.00,-4.67){};
            \node[Leaf](2)at(2.00,-4.67){};
            \node[Leaf](4)at(4.00,-4.67){};
            \node[Leaf](6)at(6.00,-4.67){};
            \node[Node](1)at(1.00,-2.33){$\GenA$};
            \node[Node](3)at(3.00,0.00){$\GenA$};
            \node[Node](5)at(5.00,-2.33){$\GenB$};
            \draw[Edge](0)--(1);
            \draw[Edge](1)--(3);
            \draw[Edge](2)--(1);
            \draw[Edge](4)--(5);
            \draw[Edge](5)--(3);
            \draw[Edge](6)--(5);
            \node[Leaf](r)at(3.00,2.33){};
            \draw[Edge](r)--(3);
        \end{tikzpicture}}
        \scalebox{.75}{
        \begin{tikzpicture}[Centering,xscale=0.17,yscale=0.24]
            \node[Leaf](0)at(0.00,-1.50){};
            \node[Leaf](2)at(2.00,-1.50){};
            \node[Node](1)at(1.00,0.00){$\GenA$};
            \draw[Edge](0)--(1);
            \draw[Edge](2)--(1);
            \node[Leaf](r)at(1.00,1.5){};
            \draw[Edge](r)--(1);
        \end{tikzpicture}}}
    + \Realization_{\Alphabet_{\Signature, \N}} \App
    \BasisE_{
        \scalebox{.75}{
        \begin{tikzpicture}[Centering,xscale=0.17,yscale=0.18]
            \node[Leaf](0)at(0.00,-4.67){};
            \node[Leaf](2)at(2.00,-4.67){};
            \node[Leaf](4)at(4.00,-4.67){};
            \node[Leaf](6)at(6.00,-4.67){};
            \node[Node](1)at(1.00,-2.33){$\GenB$};
            \node[Node](3)at(3.00,0.00){$\GenA$};
            \node[Node](5)at(5.00,-2.33){$\GenA$};
            \draw[Edge](0)--(1);
            \draw[Edge](1)--(3);
            \draw[Edge](2)--(1);
            \draw[Edge](4)--(5);
            \draw[Edge](5)--(3);
            \draw[Edge](6)--(5);
            \node[Leaf](r)at(3.00,2.33){};
            \draw[Edge](r)--(3);
        \end{tikzpicture}}
        \scalebox{.75}{
        \begin{tikzpicture}[Centering,xscale=0.17,yscale=0.24]
            \node[Leaf](0)at(0.00,-1.50){};
            \node[Leaf](2)at(2.00,-1.50){};
            \node[Node](1)at(1.00,0.00){$\GenA$};
            \draw[Edge](0)--(1);
            \draw[Edge](2)--(1);
            \node[Leaf](r)at(1.00,1.5){};
            \draw[Edge](r)--(1);
        \end{tikzpicture}}}
    + \Realization_{\Alphabet_{\Signature, \N}} \App
    \BasisE_{
        \scalebox{.75}{
        \begin{tikzpicture}[Centering,xscale=0.17,yscale=0.18]
            \node[Leaf](0)at(0.00,-4.67){};
            \node[Leaf](2)at(2.00,-4.67){};
            \node[Leaf](4)at(4.00,-4.67){};
            \node[Leaf](6)at(6.00,-4.67){};
            \node[Node](1)at(1.00,-2.33){$\GenA$};
            \node[Node](3)at(3.00,0.00){$\GenB$};
            \node[Node](5)at(5.00,-2.33){$\GenA$};
            \draw[Edge](0)--(1);
            \draw[Edge](1)--(3);
            \draw[Edge](2)--(1);
            \draw[Edge](4)--(5);
            \draw[Edge](5)--(3);
            \draw[Edge](6)--(5);
            \node[Leaf](r)at(3.00,2.33){};
            \draw[Edge](r)--(3);
        \end{tikzpicture}}
        \scalebox{.75}{
        \begin{tikzpicture}[Centering,xscale=0.17,yscale=0.24]
            \node[Leaf](0)at(0.00,-1.50){};
            \node[Leaf](2)at(2.00,-1.50){};
            \node[Node](1)at(1.00,0.00){$\GenA$};
            \draw[Edge](0)--(1);
            \draw[Edge](2)--(1);
            \node[Leaf](r)at(1.00,1.5){};
            \draw[Edge](r)--(1);
        \end{tikzpicture}}}.
    \notag
\end{align}

\section*{Open questions and future work}
We have introduced polynomial realizations of $\NaturalHopfAlgebra \App \Operad$ and used
them to establish links between this Hopf algebra and other ones. Here are open questions
raised by this work.
\begin{enumerate}[label={\bf (\arabic*)}]
    \item This concerns the further study of easterly wind lattices, specifically to
    investigate properties like semidistributivity, number of covering relations, and number
    of intervals.
    \item Another question is to delineate conditions on the operad $\Operad$ for the fact
    that $\NaturalHopfAlgebra \App \Operad$ is self-dual. The $\BasisF$ and $\BasisH$-bases
    (when $\Operad$ is free) of $\NaturalHopfAlgebra \App \Operad$, together with its
    polynomial realization could intervene to solve this question.
    \item This axis is about the construction of a Hopf algebra $\WQSym^{(s)}$ that is to
    $\WQSym$ what $\NCSF^{(s)}$ is to $\NCSF$ with the aim to obtain a stronger version of
    Theorem~\ref{thm:map_to_wqsym}.
    \item The question of understanding the structure of $\NaturalHopfAlgebra \App
    \MAs_{\Signature}$ in the general case, when $\Signature$ is any signature (having
    possibly elements of arity $0$) as this Hopf algebra generalizes the Faà du Bruno and
    the multi-symmetric function Hopf algebras.
\end{enumerate}

\printbibliography

\end{document}